\documentclass{article}
\usepackage{amsmath,amssymb}
\usepackage{graphicx}
\usepackage{float,afterpage}
\usepackage[notref,notcite]{showkeys}

\newtheorem{theorem}{Theorem}[section]
\newtheorem{rem}[theorem]{Remark}

\newtheorem{cor}[theorem]{Corollary}
\newtheorem{exmp}[theorem]{Example}
\newtheorem{definition}[theorem]{Definition}
\newtheorem{prop}[theorem]{Proposition}
\newenvironment{remark}{\begin{rem}\rm}{\end{rem}}

\renewcommand{\epsilon}{\varepsilon}
\renewcommand{\tilde}{\widetilde}
\renewcommand{\hat}{\widehat}

\newcommand{\sign}{\mathop{\mathrm{sign}}}

\usepackage{epsfig}
\title{Spectral perturbation bounds for selfadjoint
operators I\thanks{This work was partly done
during the author's stay at the University of Split, Faculty of
Electrotechnical Engineering, Mechanical Engineering and Naval Archtecture
while supported by the National Foundation
for Science, Higher Education and Technological Development of the
Republic of Croatia. Both the Foundation support and the kind hospitality
of professor Slapni\v{c}ar are gratefully acknowledged.}}
\author{Kre\v simir Veseli\'c\thanks{Fernuniversit\" at Hagen,
Fakult\" at f\" ur Mathematik und Informatik
Postfach 940, D-58084 Hagen, Germany, e-mail:
kresimir.veselic@fernuni-hagen.de.}}
\date{ }
\begin{document}
\maketitle
\begin{abstract}
We give general spectral and eigenvalue perturbation bounds for
a selfadjoint operator perturbed in the sense of the
pseudo-Friedrichs extension. We also give
several generalisations of the aforementioned extension.
The spectral bounds for finite eigenvalues are obtained by using
analyticity and monotonicity properties
(rather than variational
principles) and they are general enough to
include eigenvalues in gaps of the essential spectrum.
\end{abstract}

\setcounter{equation}{0}
\setcounter{figure}{0}
\section{Introduction}

The main purpose of this paper is to derive spectral and 
eigenvalue bounds for selfadjoint operators. If a selfadjoint operator
\(H\) in a Hilbert space \({\cal H}\) is perturbed into
\begin{equation}\label{T=H+A}
T = H + A
\end{equation}
with, say, a bounded \(A\) then the well-known spectral
spectral inclusion holds
\begin{equation}\label{inclusion}
\sigma(T) \subseteq \left\{\lambda:\
 \mbox{dist}(\lambda,\sigma(H)) \leq \|A\|\right\}.
\end{equation}
Here \(\sigma\) denotes the spectrum of a linear operator.
(Whenever not otherwise stated we shall follow the notation and 
the terminology of \cite{kato}.)

If \(H\), \(A\), \(T\) are finite Hermitian matrices then (\ref{T=H+A})
implies
\begin{equation}\label{ev_bound_normA}
|\mu_k - \lambda_k| \leq \|A\|,
\end{equation}
where \(\mu_k,\lambda_k\) are the non-increasingly ordered
eigenvalues of \(T,H\), respectively. (Here and henceforth we 
count the eigenvalues together with their multiplicities.)

Whereas (\ref{inclusion}) may be called an upper semicontinuity
bound the estimate (\ref{ev_bound_normA}) contains an
{\em existence} statement: each of the intervals
\([\lambda_k -\|A\|,\ \lambda_k +\|A\|]\) contains 'its own'
\(\mu_k\). Colloquially, bounds like (\ref{inclusion})  may be called
'one-sided' and those like (\ref{ev_bound_normA}) 'two-sided'.
As it is well-known (\ref{ev_bound_normA}) can be refined to
another two-sided bound
\begin{equation}\label{ev_bound_A+-}
\min\sigma(A) \leq \mu_k - \lambda_k \leq \max\sigma(A).
\end{equation}
In \cite{vessla1} the following 'relative' two-sided bound was derived
\begin{equation}\label{vessla1_bound}
|\mu_k - \lambda_k| \leq b|\lambda_k|,
\end{equation}
provided that
\[
|(A\psi,\psi)| \leq b(|H|\psi,\psi),\quad b < 1.
\]
This bound was found to be relevant for numerical
computations. Combining (\ref{ev_bound_normA}) and
(\ref{vessla1_bound}) we obtain
\begin{equation}\label{ev_bound2}
|\mu_k - \lambda_k| \leq a + b|\lambda_k|,
\end{equation}
or, equivalently,
\begin{equation}\label{ev_bound}
\lambda_k -a - b |\lambda_k| \leq \mu_k \leq \lambda_k + a + b|\lambda_k|,
\end{equation}
provided that
\begin{equation}\label{katobound0}
|(A\psi,\psi)| \leq a\|\psi\|^2 + b(|H|\psi,\psi),\quad b < 1.
\end{equation}

One of our goals is to extend the bound (\ref{ev_bound2}) to general
selfadjoint operators. Since these may be unbounded we have to make
precise what we mean by the sum (\ref{T=H+A}). Now, the condition
(\ref{katobound0}) is exactly the one which guarantees the existence
and the uniqueness of a closed extension \(T\) of \(H + A\),
if, say, \({\cal D}(A) \supseteq {\cal D}(|H|^{1/2})\).
The operator \(T\) is called {\em the pseudo-Friedrichs extension}
of \(H + A\), see \cite{kato}, Ch.~VI.~Th.~3.11.
Further generalisations of this construction are contained in
\cite{faris,nenciu,gkmv}. All they allow \(A\) to 
be merely a quadratic form, so (\ref{T=H+A}) is understood as
the form sum; note that the estimate (\ref{katobound0}) concerns
just forms. Particularly striking by its simplicity is the
construction made in \cite{gkmv} for the so-called quasidefinite
operators (finite matrices with this property have been studied in
\cite{vesfact}, cf.~also the references given there).  Let \(H,A\) be bounded and, in the intuitive matrix
notation,
\begin{equation}\label{quasidef}
H = 
\left[\begin{array}{rr}
H_+  &  0\\
0    & -H_-\\
\end{array}\right],\
A =
\left[\begin{array}{rr}
0   &  B\\
B^* & 0\\
\end{array}\right],\
\end{equation}
with \(H_\pm\) positive definite.
Then
\begin{equation}\label{schur}
T =
\left[\begin{array}{rr}
1            &      0     \\
B^*H_+^{-1}  &      1     \\
\end{array}\right]
\left[\begin{array}{rr}
H_+  &  0\\
0         & -H_- - B^*H_+^{-1}B\\
\end{array}\right]
\left[\begin{array}{rr}
1  &      H_+^{-1}B\\
0  &              1\\
\end{array}\right]
\end{equation}
with an obvious bounded inverse.
This is immediately transferable to unbounded 
\(H,A\) provided that
\(F = H_+^{-1/2}BH_-^{-1/2}\) is bounded. Indeed,
then (\ref{schur}) can be rewritten as
\begin{equation}\label{schur_unb}
T =
|H|^{1/2}
\left[\begin{array}{rr}
1            &      0     \\
F^*  &      1     \\
\end{array}\right]
\left[\begin{array}{rr}
1  &  0\\
0         & -1 - F^*F\\
\end{array}\right]
\left[\begin{array}{rr}
1  &      F\\
0  &              1\\
\end{array}\right]
|H|^{1/2}
\end{equation}
which is selfadjoint as a product of factors which have
bounded inverses. Note that in (\ref{katobound0})
we have \(a = 0\) and \(b = \|F\|\) and the latter need
not be less than one! 

In fact, our first task will be to
derive further constructions of operators defined as form sums.
One of them takes in (\ref{quasidef})
\[
A =
\left[\begin{array}{rr}
A_+  &  B\\
B^* &  A_-\\
\end{array}\right],
\]
where \(A_\pm\) are \(H_\pm\)-bounded as in (\ref{katobound0}).
So, we require \(b < 1\) only for 'diagonal blocks'. Another one
exhibits 'off-diagonal dominance' inasmuch as \(H_\pm\) in
(\ref{quasidef})
are a sort of \(B\)-bounded. All these constructions as well as
those from \cite{kato,faris,nenciu,gkmv} are shown to be contained
in a general abstract theorem which also helps to get a unified view
of the material scattered in the literature. This is done in Sect.~2.

As a rule each such construction will also contain a spectral
inclusion like (\ref{inclusion}). In Sect.~3 we will give some
more inclusion theorems under the condition (\ref{katobound0})
as an immediate preparation for eigenvalue estimates. 
In the proofs the quasidefinite structure will be repeatedly used.
Moreover, the decomposition (\ref{schur}) and the corresponding
invertibility property will be carried over to the Calkin algebra,
thus allowing tight control of the spectral movement including
the monotonicity in gaps both for the total and the essential spectra.

In Sect.~4 we consider two-sided bounds for finite eigenvalues. They are
obtained by using analyticity and monotonicity properties.\footnote{
Another possible approach to the monotonicity could be to use
variational principles valid also in spectral gaps,
see e.g.~\cite{griesemer} or \cite{dolbeault} but we found
the analyticity more elegant.} In order to do this
we must be able
\begin{itemize}
\item[(i)] to count the eigenvalues (note that we may be in
a gap of the essential spectrum) and
\item[(ii)] to keep the essential spectrum away from the considered
region.
\end{itemize}
The condition (i) is achieved by requiring that at least
one end of the considered interval be free from spectrum
during the perturbation (we speak od 'impenetrability').
This will be guaranteed by one of the spectral inclusion theorems
mentioned above. Similarly, (ii) is guaranteed by analogous inclusions
for the essential spectrum. Based on this we first prove
a monotonicity result for a general class of selfadjoint
holomorphic families and then establish the bound
(\ref{ev_bound2}) as well as an analogous relative bound generalising
(\ref{ev_bound_A+-}) which includes the monotonicity of eigenvalues
in spectral gaps. Another result, perhaps even more important
in practice, is the one in which the form \(A\) is perturbed
into \(B\) with \(B - A\) small with respect to \(A\)
(this corresponds to relatively small perturbations of 
the potential in quantum mechanical applications).
In this case the necessary impenetrability is obtained
by a continuation argument which assumes the knowledge of the
whole family \(H + \eta A\) instead of the mere unperturbed operator
\(H + A\). All our eigenvalue bounds are sharp.

The corresponding eigenvector bounds as well as systematic 
study of applications to various particular classes of 
operators will be treated in forthcoming papers.\\

{\small {\bf Acknowledgements.} The author is indebted to
V.~Enss, L.~Grubi\v{s}i\'{c}, R.~Hryniv, W.~Kirsch, V.~Kostrykin,
I.~Slapni\v{c}ar and I.~Veseli\'{c} for their helpful
discusions. He is also indebted to an anonymous referee
whose comments have greatly helped in the preparation 
of the final version of this paper.
}

\section{Construction of operators}
Here we
will give various constructions of selfadjoint operators by means of forms
(cf. \cite{kato,faris,nenciu,gkmv}). Sometimes our results will
generalise the aforementioned
ones only slightly, but we will still
give the proofs because their ingredients will be used
in the later work.
We shall include non-symmetric perturbations
whenever the proofs naturally allow such possibility.

\begin{definition}\label{gaps} We say that
the open interval \((\lambda_-,\lambda_+)\) is a {\em spectral gap} 
of a selfadjoint operator \(H\), if this interval belongs
to the resolvent set \(\rho(H)\) and its ends, if finite,
belong to the spectrum \(\sigma(H)\). The {\em essential spectral gap}
is defined analogously.
\end{definition}
\begin{definition}\label{representation}
We say that a sesquilinear form \(\tau\), defined
in a Hilbert space \({\cal H}\) on a dense domain \({\cal D}\)
{\em represents} an operator \(T\), if
\begin{equation}\label{represent1}
T\mbox{ is closed and densely defined,}
\end{equation}
\begin{equation}\label{represent2}
{\cal D}(T), {\cal D}(T^*) \subseteq {\cal D}
\end{equation}
\begin{equation}\label{represent3}
(T\psi,\phi) = \tau(\psi,\phi),\quad
\psi \in {\cal D}(T),\ \phi \in {\cal D},
\end{equation}
\begin{equation}\label{represent4}
(\psi,T^*\phi) = \tau(\psi,\phi),\quad
\psi \in {\cal D},\ \phi \in {\cal D}(T^*).
\end{equation}
\end{definition}
\begin{prop}\label{uniqueness}
A closed, densely defined operator \(T\) is uniquely defined
by (\ref{represent1}) -- (\ref{represent4}).
\end{prop}
{\bf Proof.} Suppose that \(T_1\) satisfies
(\ref{represent1}) -- (\ref{represent4}). Then
\[
(T\psi,\phi) = \tau(\psi,\phi) =
(\psi,T_1^*\phi),\quad \psi \in {\cal D}(T),
\ \phi \in {\cal D}(T_1^*),
\]
\[
(T_1\psi,\phi) = \tau(\psi,\phi) =
(\psi,T^*\phi),\quad \psi \in {\cal D}(T_1),
\ \phi \in {\cal D}(T^*).
\]
The first relation implies \(T_1 \supseteq T\) and the second \(T \supseteq T_1\).
Q.E.D.\\

Let \(H\) be selfadjoint in a Hilbert space \({\cal H}\) and
let \(\alpha(\cdot,\cdot)\) be a
sesquilinear form defined on \({\cal D}\) such that
\begin{equation}\label{alphabound}
|\alpha(\psi,\phi)| \leq \|H_1^{1/2}\psi\|\|H_1^{1/2}\phi\|
\quad \psi, \phi \in {\cal D}
\end{equation}
where \({\cal D}\) is a core for \(|H|^{1/2}\) and
\begin{equation}\label{H1}
H_1 = a + b|H|,\quad a,b \mbox{ real, }\
b \geq 0,\quad H_1 \mbox{ positive definite.}
\end{equation}
Then the formula
\begin{equation}\label{C}
(C\psi,\phi) = \alpha(H_1^{-1/2}\psi,H_1^{-1/2}\phi),\quad
\psi, \phi \in {\cal D},
\end{equation}
defines a \(C \in {\cal B}({\cal H})\) with
\begin{equation}\label{Cnorm}
\|C\| \leq 1
\end{equation}
(note that \(H_1^{1/2}{\cal D}\) is dense in \({\cal H}\)). The form
\(\alpha\) can obviously be extended to the form \(\alpha_{\cal Q}\),
defined on the subspace
\begin{equation}\label{Q}
{\cal Q} = {\cal D}(|H|^{1/2}) = {\cal D}(H_1^{1/2})
\end{equation}
by the formula
\begin{equation}\label{alphaQ}
\alpha_{\cal Q}(\psi,\phi) = \lim_{n,m \to \infty}
\alpha(\psi_n,\phi_n)
\end{equation}
for any sequence \(\psi_n \to \psi\), \(\phi_m \to \phi\),
\(H_1^{1/2}\psi_n \to H_1^{1/2}\psi\), \(H_1^{1/2}\phi_m \to H_1^{1/2}\phi\). Then
(\ref{alphabound}) holds for \(\alpha_{\cal Q}\) on \({\cal Q}\) and
\begin{equation}\label{CC}
(C\psi,\phi) = \alpha_{\cal Q}(H_1^{-1/2}\psi,H_1^{-1/2}\phi),\quad
\psi, \phi \in {\cal H}.
\end{equation}
The sesquilinear form for \(H\) is defined on
\({\cal Q}\) as
\begin{equation}\label{formh}
h(\psi,\phi) = (J|H|^{1/2}\psi,|H|^{1/2}\phi)
\end{equation}
with
\begin{equation}\label{Jsign}
J = \sign H.
\end{equation}
In general there may be several different sign
functions \(J\) of \(H\) with \(J^2 = 1\). The form \(h\) does not
depend on the choice of \(J\).
\begin{theorem}\label{pseudoF}
Let \(H\), \(\alpha\), \(C\), \({\cal D}\)
\({\cal Q}\) be as above and
such that
\begin{equation}\label{Czeta}
C_\zeta = (H -\zeta)H_1^{-1} + C
\end{equation}
is invertible in \({\cal B}({\cal H})\) for some \(\zeta \in {\bf C}\). Then the
form
\begin{equation}\label{formtau}
\tau = h + \alpha_{\cal Q}
\end{equation}
represents a unique closed densely defined operator \(T\) whose domain
is a core for \(|H|^{1/2}\) and which is given by
\begin{equation}\label{T-zeta}
T - \zeta = H_1^{1/2}C_\zeta H_1^{1/2},
\end{equation}
\begin{equation}\label{T*-zeta}
T^* - \overline{\zeta} = H_1^{1/2}C_\zeta^*H_1^{1/2},
\quad \zeta \in {\bf C}
\end{equation}
and, whenever \(C_\zeta^{-1} \in {\cal B}({\cal H})\),
\begin{equation}\label{invT-zeta}
(T - \zeta)^{-1} = H_1^{-1/2}C_\zeta^{-1}H_1^{-1/2} \in
{\cal B}({\cal H}),
\end{equation}
\begin{equation}\label{invT*-zeta}
(T^* - \overline{\zeta})^{-1} = H_1^{-1/2}C_\zeta^{-*}H_1^{-1/2} \in
{\cal B}({\cal H}).
\end{equation}
We call \(T\)
{\em the form sum} of \(H\) and \(\alpha\) and write
\begin{equation}\label{THalpha}
T = H + \alpha.
\end{equation}
If \(\alpha\) is symmetric then \(T\) is selfadjoint.
\end{theorem}
{\bf Proof.} In view of what was said above we may obviously suppose that
\({\cal D}\) is already equal to \({\cal Q}\).\footnote{This assumption
will be made throughout the rest of the paper,
if not stated otherwise.} We first prove that
\({\cal D}(H_1^{1/2} C_\zeta H_1^{1/2})\) is independent
of \(\zeta\) and is dense
in \({\cal H}\). Indeed, for
\(\zeta,\zeta' \in {\bf C}\) and
\(\psi \in {\cal D}(H_1^{1/2} C_\zeta H_1^{1/2})\)
we have \(\psi \in {\cal Q}\) and
\[
{\cal Q} \ni C_\zeta H_1^{1/2}\psi =
(H - \zeta)H_1^{-1}H_1^{1/2}\psi + CH_1^{1/2}\psi
\]
\[
= (H - \zeta')H_1^{-1/2}\psi + CH_1^{1/2}\psi
+ (\zeta' - \zeta)H_1^{-1/2}\psi
\]
\[
= C_{\zeta'} H_1^{1/2}\psi + (\zeta' - \zeta)H_1^{-1/2}\psi.
\]
Thus, by \((\zeta' - \zeta)H_1^{-1/2}\psi \in {\cal Q}\) we have
\(C_{\zeta'} H_1^{1/2}\psi \in {\cal Q}\), hence
\(\psi \in {\cal D}(H_1^{1/2} C_{\zeta'} H_1^{1/2})\).
Since \(\zeta,\zeta'\) are arbitrary \({\cal D}(H_1^{1/2}C_\zeta H_1^{1/2})\)
is indeed independent of \(\zeta\) and (\ref{T-zeta}) holds. Now take
\(\zeta\) with \(C_\zeta^{-1} \in {\cal B}({\cal H})\). Then the
three factors on the right hand side of (\ref{T-zeta}) have bounded, everywhere
defined inverses, so (\ref{invT-zeta}) holds as well and \(T\)
is closed.
We now prove that
\({\cal D}(T)\) is a core for \(|H|^{1/2}\) or, equivalently,
for \(H_1^{1/2}\). That is, \(H_1^{1/2}{\cal D}(T)\) must
be dense in \({\cal H}\) (see \cite{kato} III, Exercise 51.9).
By taking \(\zeta\) with
\(C_\zeta^{-1} \in {\cal B}({\cal H})\) we have
\[
H_1^{1/2}{\cal D}(T) = H_1^{1/2}{\cal D}(T - \zeta) =
\]
\[
H_1^{1/2}\left\{\psi \in {\cal Q}:\ C_\zeta H_1^{1/2}\psi
\in {\cal Q}\right\} = C_\zeta^{-1}{\cal Q}
\]
and this is dense because \(C_\zeta\) maps bicontinuously
\({\cal H}\) onto itself. In
particular, \({\cal D}(T)\) is dense in \({\cal H}\).
By
\[
C_\zeta^* = (H - \overline{\zeta})H_1^{-1} + C^*
\]
all properties derived above are seen to hold for \(T^*\) as well.
The identities (\ref{represent3}), (\ref{represent4}) follow immediately from
(\ref{T-zeta}) by using the obvious identity
\begin{equation}\label{formtauC}
\tau(\psi,\phi) - \zeta(\psi,\phi)
= (C_\zeta H_1^{1/2}\psi,H_1^{1/2}\phi),
\end{equation}
valid for any \(\psi,\phi \in {\cal Q}\),
\(\zeta \in {\bf C}\). Finally, if \(\alpha\) is symmetric then
\(T,T^*\) is also symmetric and therefore selfadjoint. Q.E.D.\\
\begin{cor}\label{pseudoFC}
Let \(H,\ \alpha,\ \tau,\ T\) be as in Theorem \ref{pseudoF}.
Then
\[
\tau(\psi,\phi) = \zeta(\psi,\phi)
\]
for some \(\psi \in {\cal Q},\ \zeta \in {\bf C}\) and
all \(\phi \in {\cal Q}\) is equivalent to
\[
\psi \in {\cal D}(T),\quad T\psi = \zeta \psi.
\]
\end{cor}
\begin{remark}\label{bounded_difference}
Although fairly general, the preceding theorem does not cover
all relevant form representations. If \(T = H + \alpha\) and
\(\alpha_1\) is any bounded form, then \(\tau_1 = \tau + \alpha_1\)
obviously generates a \(T_1\) in the sense of Def.~\ref{representation}
--- we again write \(T = H + \alpha + \alpha_1\) --- while
\(H,\ \alpha + \alpha_1\) need not satisfy the conditions
of Theorem \ref{pseudoF}.
\end{remark}
\begin{remark}\label{diag_form_bound} If \(\alpha\) is symmetric then
(\ref{alphabound}) is equivalent to
\begin{equation}\label{alphaboundsymm}
|\alpha(\psi,\psi)| \leq \|H_1\psi\|^2.
\end{equation}
In general, (\ref{alphaboundsymm}) implies (\ref{alphabound}) but with \(b\)
replaced by \(2b\).
\end{remark}
\begin{remark}\label{indep_of_ab} By Proposition
\ref{uniqueness} the operator \(T = H + \alpha\) does
not depend on the choice of \(a,b\) in the operator \(H_1\) from
(\ref{H1}). Moreover, in the construction (\ref{T-zeta})
\(H_1\) may be replaced by any selfadjoint \(H_2 = f(H)\) where
\(f\) is a real positive-valued function and
\[
0 < m \leq \frac{a + b|\lambda|}{f(\lambda)} \leq M < \infty.
\]
Then
\begin{equation}\label{invT-zeta_2}
(T - \zeta)^{-1} = H_2^{-1/2}D_\zeta^{-1}H_2^{-1/2},
\end{equation}
where
\[
D_\zeta = (H -\zeta)H_2^{-1} + D,
\]
\[
D = H_1^{1/2}f(H)^{-1/2}CH_1^{1/2}f(H)^{-1/2}
\]
and
\[
D_\zeta = H_1^{1/2}f(H)^{-1/2}C_\zeta H_1^{1/2}f(H)^{-1/2}
\]
is invertible in \({\cal B}({\cal H})\), if and only if
\(C_\zeta\) is such.
\end{remark}

%
\begin{cor}\label{pseudoF_nenciu}
Let \(H\), \(H_1 = a + |H|\), \({\cal Q}\), \(C\)
\(\alpha = \alpha_{\cal Q}\), \(h\)
and \(J\) be as in  (\ref{Q}) -- (\ref{Jsign})
such that \(J + C\)
is invertible in \({\cal B}({\cal H})\).
Then the
form \(\tau = h + \alpha\)
represents a unique closed densely defined operator
 \(T = H + \alpha\)
in the sense of Remark \ref{bounded_difference}.
Moreover, \({\cal D}(T)\) is
a core for \(|H|^{1/2}\) and
\begin{equation}\label{Tnenciu}
T + aJ = H_1^{1/2}(J + C)H_1^{1/2}
\end{equation}
(and similarly for \(T^*\)).
\end{cor}
Note that the preceding construction --- in contrast to the
related one in Theorem \ref{pseudoF} does not give an immediate
representation of the resolvent, except, if \(a = 0\).\\

In the following theorem we will use
the well known formulae
\begin{equation}\label{AB_BA}
\sigma(AB)\setminus\{0\} = \sigma(BA)\setminus\{0\},
\end{equation}
\begin{equation}\label{BA_res_AB}
(\lambda - BA)^{-1} = \frac{1}{\lambda} +
\frac{B(\lambda - AB)^{-1}A}{\lambda}
\end{equation}
\begin{equation}\label{BfAB_fBAB}
Bf(AB) = f(BA)B,
\end{equation}
where \(A,B \in {\cal B}({\cal H})\) and \(f\)
is analytic on \(\sigma(AB) \cup \{0\}\).
\begin{theorem}\label{pseudoF_nenciu_factor}
Let \(H\), \(\alpha\), \({\cal Q}\), \(C\) satisfy
(\ref{alphabound}) -- (\ref{CC}). Let, in addition,
\begin{equation}\label{Z1Z2}
C = Z_2^*Z_1,\ Z_{1,2} \in {\cal B}({\cal H}).
\end{equation}
Then \(C_\zeta\) from (\ref{Czeta}) is invertible in
\({\cal B}({\cal H})\), if and only if
\begin{equation}\label{Fzeta}
F_\zeta = 1 + Z_1H_1(H - \zeta)^{-1}Z_2^*.
\end{equation}
is such. In this case Theorem \ref{pseudoF} holds and
\begin{equation}\label{invT-zeta_Z1Z2}
(T - \zeta)^{-1} = (H - \zeta)^{-1} - 
H_1^{1/2}(H - \zeta)^{-1}Z_2^*F_\zeta^{-1}Z_1H_1^{1/2}(H - \zeta)^{-1}.
\end{equation}

\end{theorem}
{\bf Proof.} \(C_\zeta\) is invertible in
\({\cal B}({\cal H})\), if and only if
\[
1 + H_1(H - \zeta)^{-1}C = 1 + H_1(H - \zeta)^{-1}Z_2^*Z_1
\]
is invertible in \({\cal B}({\cal H})\).
Now,
\[
\sigma(H_1(H - \zeta)^{-1}Z_2^*Z_1)\setminus\{0\} =
\sigma(Z_1H_1(H - t)^{-1}Z_2^*)\setminus\{0\}
\]
Hence \(F_\zeta\) is invertible in \({\cal B}({\cal H})\)
if an only if \(C_\zeta\) is such.
In this case (\ref{invT-zeta}) gives
\[
(T - \zeta)^{-1} = H_1^{-1/2}(1 +
H_1(H - \zeta)^{-1}Z_2^*Z_1)^{-1}H_1^{1/2}(H - \zeta)^{-1} =
\]
\[
H_1^{-1/2}\left(1 - H_1(H - \zeta)^{-1}Z_2^*Z_1(1 +
H_1(H - \zeta)^{-1}Z_2^*Z_1)^{-1}\right)
H_1^{1/2}(H - \zeta)^{-1} =
\]
\[
(H - \zeta)^{-1} - H_1^{1/2}(H - \zeta)^{-1}Z_2^*F_\zeta^{-1}
Z_1H_1^{1/2}(H - \zeta)^{-1}.
\]
 Q.E.D.\\

We now apply Theorem \ref{pseudoF} to further cases in which
the key operator \(C_\zeta\) from (\ref{Czeta}) is invertible in
\({\cal B}({\cal H})\).

\begin{theorem}\label{pseudoF1}
Let \(H\) be selfadjoint and let \(\alpha\) satisfy (\ref{alphabound})
with \(b < 1\)
and (\ref{H1}). Then the conditions of Theorem \ref{pseudoF}
are satisfied and \(\zeta = \lambda +
i\eta \in \rho(T)\) whenever
\begin{equation}\label{zetaregion}
|\eta| > \frac{a + |\lambda|b.}{\sqrt{1 - b^2}}
\end{equation}
\end{theorem}
{\bf Proof.} To prove \(C_\zeta^{-1} = \left((H - \zeta)H_1^{-1} + C\right)^{-1}
\in {\cal B}({\cal H})\) it is enough to find a \(\zeta
= \lambda + i\eta\) such that
\begin{equation}\label{etaneumann}
\|H_1(H - \zeta)^{-1}\| < 1.
\end{equation}
Now,
\[
\|(H - \zeta)^{-1}H_1\| \leq \sup_{\xi \in {\bf R}}
\psi(\xi,a,b,\lambda,\eta),\quad \psi(\xi,a,b,\lambda,\eta) =
\frac{b|\xi| + a}{\sqrt{(\xi - \lambda)^2 + \eta^2}}
\]
A straightforward, if a bit tedious, calculation
(see Appendix) shows
\begin{equation}\label{maxf}
\max_{\xi}\psi =
\frac{1}{|\eta|}\sqrt{(a + |\lambda|b)^2 + b^2\eta^2}
\end{equation}
Hence (\ref{zetaregion}) implies (\ref{etaneumann}). Q.E.D.\\

Another similar criterion for the validity of
Theorem \ref{pseudoF} --- oft independent of that of
Theorem \ref{pseudoF1} is given by the following
\begin{cor}\label{V1V2cor} Let \(H\), \(\alpha\), \(C\)
satisfy (\ref{alphabound})--(\ref{C}) and let
\begin{equation}\label{Z1Z2<1}
\|Z_1H_1(H - \zeta)^{-1}Z_2^*\| < 1.
\end{equation}
for some \(\zeta \in \rho(H)\) and with \(Z_{1,2}\)
from (\ref{Z1Z2}). Then Theorem
\ref{pseudoF_nenciu_factor} applies.
\end{cor}

Typically we will have
\begin{equation}\label{V1V2}
\alpha(\psi,\phi) = (V_1\psi,V_2\phi),
\end{equation}
where \(V_{1,2}\) are linear operators defined on \({\cal Q}\)
such that
\begin{equation}\label{Z1Z2H1}
Z_{1,2} = V_{1,2}H_1^{1/2} \in {\cal B}({\cal H}).
\end{equation}
In this case the formula (\ref{Z1Z2<1}) can be given a
more familiar, if not always rigorous, form (cf.~\cite{nenciu})
\[
\|V_1(H - \zeta)^{-1}V_2^*\| < 1.
\]
\begin{remark}\label{mygeneralisations}
If \(\alpha(\psi,\phi) = (A\psi,\phi)\), where \(A\) is a linear operator
defined on \({\cal D} \subseteq {\cal D}(H)\), \({\cal D}\) a core for
\(|H|^{1/2}\) then Theorem \ref{pseudoF1}
applies and, by construction, the obtained operator coincides
with the one in \cite{kato} VI.~Th.~3.11. The uniqueness of \(T\)
as an extension of \(H + A\), proved in \cite{kato} makes no sense
in our more general, situation.
Our notion of form uniqueness (which was used by
\cite{nenciu} in the symmetric case) will be appropriate
in applications to both Quantum and Continuum Mechanics. Thus, our Theorem
\ref{pseudoF1} can be seen as a slight generalisation of
\cite{kato}. 
On the other side, our proof of Theorem 
\ref{pseudoF} closely follows the one from \cite{kato}.

Cor.~\ref{pseudoF_nenciu} and  Th.~\ref{pseudoF_nenciu_factor}
are essentially Theorems.~2.1, 2.2 in \cite{nenciu}
except for the following: (i)
our \(\alpha\) need not be symmetric, (ii) we use 
a more general factorisability (\ref{Z1Z2}) instead of (\ref{V1V2})
which is supposed in \cite{nenciu} and finally, (iii)
we need no relative compactness argument to establish
Theorem \ref{pseudoF_nenciu_factor}.
The fact that the mentioned results from \cite{nenciu}
are covered by our theory will facilitate to handle
perturbations of the form \(\alpha\) which are not easily
accessible, if \(\alpha\) is factorised as in
(\ref{V1V2}). The spectral inclusion
formula (\ref{zetaregion}) seems to be new.

Thus, our Theorem \ref{pseudoF} seems to cover essentially all
known constructions thus far.\footnote{There are two obvious
extensions: (i) adding a bounded form (Remark \ref{bounded_difference})
and (ii) multiplying \(T\) by a bicontinuous operator. An example of the latter is \(T = H + \alpha\) described in Cor.~\ref{pseudoF_nenciu}.}
\end{remark}
Next we give some results on the invariance of the essential
spectrum.
\begin{theorem}\label{relcomp}
Let \(H\), \(h\), \(\alpha\), \(C\), \({\cal D}\)
\({\cal Q}\) satisfy (\ref{alphabound}) -- (\ref{Jsign}) with
\(\alpha\) symmetric.

(i) If the operator \(C\) is compact then Theorem
\ref{pseudoF1} holds and \(\sigma_{ess}(T) =
\sigma_{ess}(H)\). (ii) If Theorem \ref{pseudoF}
holds and \(H_1^{-1}C\) is compact then again
\(\sigma_{ess}(T) = \sigma_{ess}(H)\).
\end{theorem}
{\bf Proof.}
In any of the cases (i), (ii) we can find a \(\zeta\) for which
\(C_\zeta^{-1} \in {\cal B}({\cal H})\) (in the case (i) this follows
from the known argument that for a compact \(C\) the estimate
(\ref{alphabound})
will hold with arbitrarily small \(b\)) so Theorem \ref{pseudoF}
holds anyway. By (\ref{invT-zeta}) we have
\[
(T - \zeta)^{-1} - (H - \zeta)^{-1} =
\]
\[
H_1^{-1/2}\left(((H - \zeta)H_1^{-1} + C)^{-1} -
H_1(H - \zeta)^{-1} \right)H_1^{-1/2} =
\]
\[
H_1^{-1/2}\left((1 + A)^{-1} - 1 \right)H_1^{1/2}(H - \zeta)^{-1}
\]
where \(H_1^{-1}A = (H - \zeta)^{-1}C\) is compact and by
\(C_\zeta^{-1} \in {\cal B}({\cal H})\) also
\((1 + A)^{-1} \in {\cal B}({\cal H})\). Hence
\[
(T - \zeta)^{-1} - (H - \zeta)^{-1} =
-H_1^{-1/2}A(1 + A)^{-1}H_1^{1/2}(H - \zeta)^{-1}
\]
is compact and the Weyl theorem applies. Q.E.D.\\

Finally we borrow from \cite{nenciu} the following result
which will be
of interest for
Dirac operators with strong Coulomb potentials.

\begin{theorem}\label{relcomp_nenciu_factor} Let
\(H\), \(\alpha\), \({\cal Q}\), \(C\), \(V_{1,2}\),
\(Z_{1,2}\), \(T\) be as in Theorem \ref{pseudoF_nenciu_factor}.
Let, in addition, \(C_\zeta\) from Theorem \ref{pseudoF}
be invertible in \({\cal B}({\cal H})\) and
\begin{enumerate}
\item \(H\) have a bounded inverse,
\item the operator \(Z_2^*(H - \zeta)^{-1}Z_1\)
be compact for some (and then all) \(\zeta \in \rho(H)\).
\end{enumerate}
Then 
\(\sigma_{ess}(T) \subseteq \sigma_{ess}(H)\).
\end{theorem}

The key invertibility of the operator \(C_\zeta\) can be achieved
in replacing the requirement \(b < 1\) in (\ref{alphabound}) by some
condition on the structure of the perturbation. One such structure
is given, at least symbolically, by the matrix
\begin{equation}\label{matrix}
\left[\begin{array}{rr}
W_+  &  B\\
B    & -W_-\\
\end{array}\right],\
\end{equation}
where \(W_\pm\) are accretive. Such operator
matrices appear in various applications
(Stokes operator, Dirac operator, especially on a manifold
(\cite{gkmv}, \cite{winklmeier}) and the like). Even more general
cases could be of interest, namely
those where \(b < 1\) in (\ref{alphabound}) is required to hold
only on the ``diagonal blocks'' of the perturbation \(\alpha\). We have
\begin{theorem} \label{diagonalb} Let \(H\), \(\alpha =
\alpha_{\cal Q}\), \(C\), \(h\) satisfy 
(\ref{alphabound}) -- (\ref{Jsign})
such that \(H\) has a spectral gap \((\lambda_-,\lambda_+)\) containing zero.
Suppose
\begin{equation}\label{alpha+-0}
\pm \Re\alpha(\psi,\psi) \leq 
a_\pm\|\psi\|^2 + b_\pm\||H - d|^{1/2}\psi\|^2,
\ \psi \in P_\pm{\cal Q},
\end{equation}
\begin{equation}\label{ab+-}
a_\pm > 0,\ 0 < b_\pm < 1,
\end{equation}
\begin{equation}\label{alphacorner}
\alpha(\psi,\phi) = \overline{\alpha(\phi,\psi)},\quad
\psi \in P_+{\cal Q}, \phi \in P_-{\cal Q}.
\end{equation}
where \(P_\pm = (1 \pm J)/2\). Finally, suppose
\begin{equation}\label{lhatless}
\hat{\lambda}_- = \lambda_- + b_-|\lambda_-|
<
\hat{\lambda}_+ = \lambda_+ - b_+|\lambda_+|.
\end{equation}
Then \(\tau = h + \alpha\)
generates a closed, densely defined operator \(T\)
with \({\cal D}(T)\) a core for \(|H|^{1/2}\) and
\begin{equation}\label{inclusion2}
(\hat{\lambda}_-,\hat{\lambda}_+) +
i{\bf R} \subseteq \rho(T).
\end{equation}
The operator \(T\) is selfadjoint, if \(\alpha\) is symmetric.
\end{theorem}
{\bf Proof.}  We split the perturbation \(\alpha\)
into two parts
\[
\alpha = \chi + \chi'
\]
where \(\chi\) is the 'symmetric diagonal part' of \(\alpha\), that is,
\[
\alpha_d(\psi,\phi) = \alpha(P_+\psi,P_+\phi) + \alpha(P_-\psi,P_-\phi),
\]
\[
\chi(\psi,\phi) = \frac{1}{2}\left(\alpha_d(\psi,\phi)
+ \overline{\alpha_d(\phi,\psi)}\right).
\]
Symbolically,\footnote{
Throughout this paper we will freely use matrix notation for
bounded operators as well as for unbounded ones or forms whenever the latter
are unambigously defined. The matrix partition refers to the
orthogonal decomposition \({\cal H} = P_+{\cal H}\oplus P_- {\cal H}\).}
\[
\chi =
\left[\begin{array}{rr}
\chi_+   &  0 \\
0        & -\chi_- \\
\end{array}\right],
\
h =
\left[\begin{array}{rr}
h_+   &  0 \\
0        & -h_- \\
\end{array}\right].
\]
Now (\ref{alpha+-0}) and the standard perturbation result for
closed symmetric forms (\cite{kato} Ch.~VI, Th. 3.6) implies
\(\tilde{h}_\pm = h_\pm + \chi_\pm\) is symmetric,
bounded from below by
\[
\pm \lambda_\pm - b_\pm|\lambda_\pm| - a_\pm
\]
and closed on \({\cal Q}\). The thus generated selfadjoint
operator \(\tilde{H}_\pm \) has 
\({\cal D}(|\tilde{H}_\pm |^{1/2}) = P_\pm{\cal Q}\).
Now,
\[
\tau = h + \alpha = \tilde{h} + \chi',\quad
\tilde{h} =
\left[\begin{array}{rr}
\tilde{h}_+   &  0 \\
0        &  -\tilde{h}_- \\
\end{array}\right].
\]
We write
\[
\tau = h + \alpha = \tilde{h} + \chi',\
\tilde{h} =
\left[\begin{array}{rr}
\tilde{h}_+   &  0 \\
0        & -\tilde{h}_- \\
\end{array}\right],\
\tilde{H} =
\left[\begin{array}{rr}
\tilde{H}_+   &  0 \\
0        & -\tilde{H}_- \\
\end{array}\right],
\]
where \(\tilde{H}\) has a spectral gap contained in
\((\tilde{\lambda}_-,\tilde{\lambda}_+)\) and
\[
J = \sign(\tilde{H} - d) = \sign H,\quad
\tilde{\lambda}_- < d < \tilde{\lambda}_+.
\]
We will apply Theorem \ref{pseudoF}
to \(\tilde{H}\), \(\chi'\).
We have first to prove that \(\tilde{H}, \chi'\) satisfy the conditions
 (\ref{alphabound}), (\ref{H1}) (possibly with different constants \(a,b\)).
By (\ref{alpha+-0}) we have
\[
0 \leq h_\pm \leq \frac{a_\pm}{1 - b_\pm} +
\frac{\tilde{h}_\pm}{1 - b_\pm}.
\]
Hence
\[
|H| \leq c|\tilde{H} - d|
\]
for any \( d \in (\tilde{\lambda}_-, \tilde{\lambda}_+)\)
and some \(c = c(d)\). So, \(\tilde{H}, \chi'\) satisfy
(\ref{alphabound}), (\ref{H1}) with \(|H|\) replaced
by \(|\tilde{H} - d|\).
We take \(\zeta = d + i\eta\) and set
\begin{equation}\label{Ttilde}
\tilde{T} - \zeta = |\tilde{H} - d|^{1/2}
D_\zeta|\tilde{H} - d|^{1/2}
\end{equation}
\[
D_\zeta = J - \zeta|\tilde{H} - d|^{-1} + D,
\]
\[
(D\psi,\phi) =
\chi'(|\tilde{H} - d|^{-1/2}\psi,|\tilde{H} - d|^{-1/2}\phi),
\]
\[
D =
\left[\begin{array}{rr}
D_+   &     F    \\
F^*     &  -D_-  \\
\end{array}\right].
\]
By the construction we have
\begin{equation}\label{+-chi'accretive}
\Re\chi'(P_\pm\psi,P_\pm\psi) = 0.
\end{equation}
Hence 
\(D_\pm\) are skew Hermitian
and
\[
D_\zeta
= \left[\begin{array}{rr}
1 - i\eta(\tilde{H}_+ -d)^{-1} + D_+   &   F   \\
F^*    &  -1 - i(d - \tilde{H}_-)^{-1} - D_-\\
\end{array}\right]
\]
where the first diagonal block is uniformly accretive
and the second uniformly dissipative, so 
\(D_\zeta^{-1} \in {\cal B}({\cal H})\) by virtue of the
factorisation (\ref{schur}) which obviously holds in this
case, too.
Thus,
\[
(T - \zeta)^{-1} = |\tilde{H} - d|^{-1/2}
\tilde{D}_\zeta^{-1}|\tilde{H} - d|^{1/2}
\in {\cal B}({\cal H})
\]
and Theorem \ref{pseudoF} applies. Note also that \(|H|^{1/2}\)
and \(|\tilde{H}|^{1/2}\) have the same set of cores.
Q.E.D.\\
\begin{cor}\label{lift_gap}
If in the preceding theorem we drop the
condition (\ref{lhatless}) or even the existence
of the spectral gap of \(H\) we still have
\(T = H + \alpha\) but without the spectral inclusion
(\ref{inclusion2}).
\end{cor}
{\bf Proof.} We first apply the preceding theorem
to \(\hat{T} = \hat{H} + \alpha\) with \(\hat{H} = H + \delta J\)
and \(\delta > 0\) large enough to insure that
(\ref{lhatless}) holds. Then set 
\(T = \hat{T} - \delta J\). Q.E.D.

\begin{remark}\label{elegant}
Theorem \ref{diagonalb} becomes particularly elegant,
if we set \(a_\pm,b_\pm = 0\). If, in addition, \(\alpha\)
is taken as symmetric then we have a 'quasidefinite form'
\(\tau\) as was mentioned in Sect.~1. In this case
we require only the condition (\ref{alphabound}) with no
restriction on the size of \(a,b\) (for \(H,\alpha\) non-negative
this is a well-known fact).
\end{remark}

There is an
alternative proof of Theorem \ref{diagonalb} 
which we now illustrate (we assume for 
simplicity that \(a_\pm = 0\)). Instead of the pair \(H,\ \alpha\)
consider \(JH = |H|,\ J\alpha\) where \(J\alpha\) is the 'product
form' naturally defined by
\[
J\alpha(\psi,\phi) = \alpha(\psi,J\phi)
\]
As one immediately sees the new form
\[
J\tau = Jh + J\alpha
\]
is sectorial and its symmetric part \(Jh\) is closed non-negative,
so by the standard theory (\cite{kato} Ch.~VI.~\S 3) \(J\tau\)
generates a closed sectorial operator which we denote by
\(JT\). Symbolically,
\[
JT =
\left[\begin{array}{rr}
1  &  0 \\
0  & -1 \\
\end{array}\right]
\left[\begin{array}{rr}
A_+  &  B\\
B^*& -A_-\\
\end{array}\right]
=
\left[\begin{array}{rr}
A_+  &  B\\
-B^*& A_-\\
\end{array}\right].
\]
The reason why we still stick
at our previous proof is its constructivity (here we have no
direct access to the resolvent) as well as its 'symmetry',
(here even for a symmetric \(\alpha\) a detour through
non-symmetric objects is made).\\


Another case in which Theorem \ref{pseudoF} can be applied is the one in which (\ref{matrix}) is 'off-diagonally dominant' (cf.~\cite{winklmeier}). We set
\begin{equation}\label{Hoff}
H =
\left[\begin{array}{rr}
0           &   B \\
B^*         &   0 \\
\end{array}\right]
\end{equation}
where \(B\) is a closed, densely defined operator between
the Hilbert spaces \({\cal H}_-\) and \({\cal H}_+\). It
is easy to see that \(H\) is selfadjoint on
\({\cal D}(B^*) \oplus {\cal D}(B)\) (see \cite{thaller}, Lemma 5.3).
Denote by
\begin{equation}\label{polarB}
B = U\sqrt{B^*B}
\end{equation}
the corresponding polar decomposition (see \cite{kato},
Ch.~VI, \S 2.7) and suppose that \(U\) is an isometry
from \({\cal H}_-\) onto \({\cal H}_+\). Then
\[
H =
\left[\begin{array}{rr}
0           &   U\sqrt{B^*B} \\
U^*\sqrt{BB^*}         &   0 \\
\end{array}\right]
=
\left[\begin{array}{rr}
0           &   U \\
U^*       &   0   \\
\end{array}\right]
\left[\begin{array}{rr}
\sqrt{B^*B}  &                \\
0            &    \sqrt{BB^*} \\
\end{array}\right]
= J|H|,
\]
\[
J^2 = I.
\]
The form \(\alpha\) is defined as follows. Denoting
\[
\psi =
\left[\begin{array}{r}
\psi_+ \\
\psi_- \\
\end{array}\right],\quad
\psi_+\in {\cal D}(B^*),\ \psi_-\in {\cal D}(B)
\]
we set
\begin{equation}\label{alpha+-}
\alpha(\psi,\phi) = \alpha_+(\psi_+,\phi_+) - \alpha_-(\psi_-,\phi_-)
\end{equation}
where \(\alpha_\pm\), defined on \({\cal D}(B^*)\), \({\cal D}(B)\),
respectively, are symmetric and non-negative.
\begin{theorem} \label{offdiag}
 Let \(H\), \(\alpha\), \(B\), \(U\)
be as above.
Let
\begin{equation}\label{alpha+bounds}
\alpha_+(\psi,\psi) \leq a_+\|\psi\|^2 + 
b_+\left((BB^*)^{1/2}\psi,\psi\right),\quad
\psi \in {\cal D}(B^*),
\end{equation}
\begin{equation}\label{alpha-bounds}
\alpha_-(\psi,\psi) \leq a_-\|\psi\|^2 +
 b_-\left((B^*B)^{1/2}\psi,\psi\right),\quad
\psi \in {\cal D}(B),
\end{equation}
for some \(a_\pm, b_\pm > 0\). Then \(\tau = h + \alpha\) generates
a unique \(T = H + \alpha\) in the sense of
Definition \ref{representation}.
\end{theorem}
{\bf Proof.} Since \(\alpha\) is defined on
\[
{\cal D}(B^*) \oplus {\cal D}(B) =
{\cal D}(\sqrt{BB^*}) \oplus {\cal D}(\sqrt{B^*B})
\]
which is obviously a core for \(|H|^{1/2}\) we can use
(\ref{alphaQ}) to extend \(\alpha\) to \({\cal Q}\)
still keeping the estimates (\ref{alpha+bounds}),
(\ref{alpha-bounds}) and similarly with \(\alpha_\pm\).
(For simplicity we denote the extended
forms again by \(\alpha,\alpha_\pm\), respectively.)

We first consider the special case, in which \(B\)
has an inverse in \({\cal B}({\cal H}_+, {\cal H}_-)\).
Then we can obviously assume that \(a_\pm = 0\) (by increasing the
size of \(b_\pm\), if necessary, note that now both \(BB^*\),
\(B^*B\) are positive definite).
Clearly, \(H^{-1} \in
{\cal B}({\cal H})\), so we may use the representation
\(T = |H|^{1/2}(J + D)|H|^{1/2}\) with
\[
J = \sign H =
\left[\begin{array}{ll}
0     &   U   \\
U^*   &   0   \\
\end{array}\right]
\]
and
\[
D =
\left[\begin{array}{ll}
D_+      &   0   \\
0        &  -D_- \\
\end{array}\right]
\]
where \(D_\pm\) are bounded symmetric non-negative. Now,
we have to prove the bounded invertibility of
\begin{equation}\label{winkl}
J + D =
\left[\begin{array}{ll}
1        &   0 \\
-D_-U^*  &   1 \\
\end{array}\right]
\left[\begin{array}{ll}
U      &   0             \\
0     &  U^* + D_-U^*D_+ \\
\end{array}\right]
\left[\begin{array}{ll}
U^*D_+      &   1   \\
1           &   0   \\
\end{array}\right]
\end{equation}
which, in turn, depends on the bounded invertibility of
\[
U^* + D_-U^*D_+
\]
or, equivalently, of \(1 + UD_-U^*D_+\).
The latter is true because
the spectrum of the product of two bounded symmetric non-negative
operators is known to be real and non-negative.

In general we first apply Theorem \ref{offdiag} to \(\tau_1 = h_1 + \alpha\)
where \(h_1\) belongs to
\[
\left[\begin{array}{rr}
0           &   B_1 \\
B_1^*         &   0 \\
\end{array}\right]
\]
and
\[
B_1 = (\sqrt{BB^*} + \delta)U,\ \delta > 0.
\]
Indeed, by \(B_1^*B_1 = B^*B + \delta\) and \(B_1B_1^* =
UB_1B_1^*U^* + \delta = BB^* + \delta\) (here we have used the assumed
isomorphy property of \(U\))
the inequalities (\ref{alpha+-})
are valid for \(B_1\) as well.
Thus, \(\tau_1\) generates \(T_1\) and \(\tau = \tau_1 + (\tau - \tau_1)\)
generates \(T\) the difference \(\tau - \tau_1\)
being bounded. Q.E.D.\\

\begin{remark}\label{offdiag_rem}
If in the preceding theorem we have \(B^{-1}
\in {\cal B}({\cal H}_+,{\cal H}_-)\)
then we can take \(a_\pm = 0\) and \(T^{-1} \in 
{\cal B}({\cal H})\) follows. This is immediately
seen from the factorisation (\ref{winkl}).
\end{remark}
\begin{remark} The property of off-diagonal dominance was
used in \cite{winklmeier} for a special
Dirac operator with a bounded form \(\alpha\) including
the decomposition (\ref{winkl}). This decomposition
has a similar disadvantage as the one described in 
Remark \ref{elegant}: it is not symmetric i.e.~it has not
the form of a congruence like e.g.~(\ref{schur}), but we
know of no better as yet.
\end{remark}
If in (\ref{matrix}) we drop the positive definiteness of, say,
\(H_-\) we still may have a positive definite Schur complement.
This gives one more possibility of constructing selfadjoint operators.
\begin{theorem} \label{schur_construction}
Let \(\tau\) be a symmetric sesquilinear form defined on a dense subspace
\({\cal Q} \subseteq {\cal H}\). Let \(P_+,P_-\) be an orthogonal decomposition
of the identity such that
\begin{itemize}
\item[(i)] \(P_\pm {\cal Q} \subseteq {\cal Q}\),
\item[(ii)] \(\tau\), restricted to \(P_+ {\cal Q}\) is closed
and positive definite,
\item[(iii)]
\[
\sup_{\psi \in P_-{\cal Q},\phi \in P_+{\cal H},\ \psi, \phi \neq 0}
\frac{|\tau(\psi,H_+^{-1/2}\phi)|}{\|\psi\|\|\phi\|}
< \infty,
\]
where \(H_+\) is the operator generated by \(\tau\) in
\(P_+{\cal H}\),
\item[(iv)] denoting by \(N \in {\cal B}(P_-{\cal H},P_+{\cal H})\)
the operator, defined by \((N\psi,\phi) = \tau(\psi,H_+^{-1/2}\phi)\),
the form
\begin{equation}\label{tauN}
P_-{\cal Q} \ni \psi,\phi \mapsto - \tau(\psi,\phi) + (N\psi,N\phi)
\end{equation}
is closed and positive definite.
\end{itemize}
Then there exists a unique selfadjoint operator \(T\) such that
\begin{itemize}
\item[(a)] \({\cal D}(T) \subseteq {\cal Q}\),
\item[(b)] \(\tau(\psi,\phi) = (T\psi,\phi),\quad \psi \in {\cal D}(T),
\phi \in {\cal Q}.\)
\end{itemize}
The operator \(T\) is given by the formulae
\begin{equation}\label{Tschur}
T = WH_1W^*,
\end{equation}
\begin{equation}\label{Wschur}
W =
\left[\begin{array}{cc}
1            &     0 \\
NH_+^{-1/2}  &     1 \\
\end{array}\right]
\in {\cal B}({\cal H}),
\end{equation}
\begin{equation}\label{T1schur}
H_1 =
\left[\begin{array}{cc}
H_+     &     0 \\
0       &   -\tilde{H}_- \\
\end{array}\right],
\end{equation}
where \(\tilde{H}_-\) is generated by the form (\ref{tauN}).
\end{theorem}
{\bf Proof.}  Obviously
\[
T^{-1} = W^{-*}H_1^{-1}W^{-1} \in {\cal B}({\cal H}),
\]
where every factor is bounded. Also
\[
W^*{\cal Q} \subseteq {\cal Q},\quad W^{-*}{\cal Q} \subseteq {\cal Q},
\]
\[
{\cal D}(T) \subseteq {\cal Q} = {\cal D}(|T_1|^{1/2}).
\]
Now take
\[
\psi =
\left[\begin{array}{r}
\psi_1  \\
\psi_2  \\
\end{array}\right]
\in {\cal D}(T),\quad
\left[\begin{array}{r}
\phi_1  \\
\phi_2  \\
\end{array}\right]
\in {\cal Q}.
\]
Then
\[
(T\psi,\phi) = (H_1W^*\psi,W^*\phi) =
\]
\[
\left(H_1\left[\begin{array}{r}
\psi_1 + H_+^{-1/2}N\psi_2 \\
\psi_2  \\
\end{array}\right],
\left[\begin{array}{r}
\phi_1 + H_+^{-1/2}N\phi_2 \\
\phi_2  \\
\end{array}\right]\right) =
\]
\[
\left(\left[\begin{array}{r}
H_+^{1/2}\psi_1 + N\psi_2 \\
-\tilde{H}_-^{1/2}\psi_2  \\
\end{array}\right],
\left[\begin{array}{r}
H_+^{1/2}\phi_1 + N\phi_2 \\
\tilde{H}_-^{1/2}\phi_2  \\
\end{array}\right]\right) =
\]
\[
\tau(\psi_1,\phi_1) + (H_+^{1/2}\psi_1,N\phi_2) +
(N\psi_2,H_+^{1/2}\phi_1) +
\]
\[
(N\psi_2,N\phi_2) + \tau(\psi_2,\phi_2) - (N\psi_2,N\phi_2).
\]
Now by \((N\psi_2,H_+^{1/2}\phi_1) = \tau(\psi_2,\psi_1)\)
we obtain
\[
(T\psi,\phi) = \tau(\psi,\phi)
\]
whereas the uniqueness follows from Proposition \ref{uniqueness}. Q.E.D.\\

\section{More spectral inclusions}
Some spectral inclusion results are already contained
in the construction Theorems \ref{pseudoF1} and \ref{diagonalb}.
They control the spectral gap at zero.
In the sequel we produce additional results valid
for general spectral gaps. We restrict ourselves here and in
the following to symmetric forms \(\alpha\)
and therefore to selfadjoint operators \(T = H + \alpha\).
\begin{theorem} \label{window0}
Let \((\lambda_-,\lambda_+)\) be an open interval,
contained in \(\rho(H)\) such that \(\lambda_\pm \in \sigma(H)\)
(we allow \(\lambda_\pm = \pm\infty\)) and let
\(T = H + \alpha\) satisfy Theorem \ref{pseudoF1}.
Let, in addition, the open interval
\begin{equation}\label{I0}
{\cal I} = (\lambda_- + (a + b|\lambda_-|),\ \lambda_+ - (a + b|\lambda_+|))
\end{equation}
be non-void. Then \({\cal I} \in \rho(T)\).
\end{theorem}
{\bf Proof.} Without loss
of generality we may take \(\lambda_+ > 0\) (otherwise consider
\(-H,-T\)). We supose first that both \(\lambda_-\) and \(\lambda_+\)
are finite. For \(d \in (\lambda_-, \lambda_+)\) we will have
\[
(T - d)^{-1} = H_1^{-1/2}((H - d)H_1^{-1} +  C)^{-1}H^{-1/2}
\in {\cal B}({\cal H}),
\]
if
\[
\|(H - d)^{-1}H_1\| < 1.
\]
Now,
\[
\|(H - d)^{-1}H_1\| = \sup_{\lambda \not\in (\lambda_-,\lambda_+)}f(\lambda),
\]
\[
f(\lambda) = \frac{b|\lambda| + a}{|\lambda - d|}.
\]
We now compute the supremum above.
\begin{itemize}
\item[Case 1: \(\lambda_- > 0\).] Then \(d > 0\).
\begin{equation}\label{prime1}
\lambda \geq \lambda_+:¸\quad \left(\frac{b\lambda + a}{\lambda - d}\right)'
= \frac{b(\lambda - d) - (b\lambda + a) }{(\lambda - d)^2}
= \frac{-db - a}{(\lambda - d)^2},
\end{equation}
\[
\max_{\lambda \geq \lambda_+}f(\lambda) =
 \frac{b\lambda_+ + a}{\lambda_+ - d}
> b;
\]
\begin{equation}\label{prime2}
0 \leq \lambda \leq \lambda_-:¸\quad
\left(\frac{b\lambda + a}{d - \lambda}\right)'
= \frac{b(d - \lambda) + (b\lambda + a) }{(\lambda - d)^2}
= \frac{db + a}{(\lambda - d)^2},
\end{equation}
\[
\max_{0 \leq \lambda \leq \lambda_-}f(\lambda)
= \frac{b\lambda_- + a}{d - \lambda_-}
> \frac{a}{d};
\]
\begin{equation}\label{prime3}
\lambda \leq 0:¸\quad \left(\frac{-b\lambda + a}{d - \lambda}\right)'
= \frac{-b(d - \lambda) + (-b\lambda + a) }{(\lambda - d)^2}
= \frac{-db + a}{(\lambda - d)^2},
\end{equation}
\[
\max_{\lambda \leq 0}f(\lambda) =
\left\{\begin{array}{cc}
a/d,   & a > db\\
b,     & a \leq db\\
       \end{array}\right.
\]
Altogether,
\[
\max_{\lambda \not\in (\lambda_-,\lambda_+)}f(\lambda) =
\max\left\{\frac{b\lambda_+ + a}{\lambda_+ - d},\
\frac{b\lambda_- + a}{d - \lambda_-}\right\}
\]
and this is obviously less than one, if \(d \in {\cal I}\).
\item[Case 2: \(\lambda_- \leq 0\).] Then \(d\)
may be negative. By (\ref{prime1}),
\begin{equation}\label{sup+}
\sup_{\lambda \geq \lambda_+}f(\lambda) =
\left\{\begin{array}{cc}
\frac{b\lambda_+ + a}{\lambda_+ - d}   & a + db \geq 0\\
b,     & a + db \leq 0\\
       \end{array}\right.
\end{equation}

By (\ref{prime3}),
\begin{equation}\label{sup-}
\sup_{\lambda \leq \lambda_-}f(\lambda) =
\left\{\begin{array}{cc}
\frac{-b\lambda_- + a}{d - \lambda_-}   & a > db\\
b,     & a \leq db\\
       \end{array}\right.
\end{equation}
\end{itemize}
Again, both suprema are less than one, if \(d \in {\cal I}\).
If one of \(\lambda_\pm\) is infinite the proof
goes along the same lines and is simpler still. Q.E.D.\\

Tighter bounds can be obtained, if more is known on the
perturbation \(\alpha\).
If \(\alpha\) is, say, non-negative then
\[
\alpha = \alpha_0 + e_0,\quad
e_0 = \inf_\psi\frac{\alpha(\psi,\psi)}{(\psi,\psi)}
\]
and both \(\alpha_0\) and \(e_0\) are non-negative. Now for
\[
T = H + \alpha =
H_1^{1/2}(HH_1^{-1} +  C)^{-1}H_1^{1/2}
\]
we have
\[
T - e_0 = H + \alpha_0 =
H_1^{1/2}(HH_1^{-1} +  C_0)^{-1}H_1^{1/2}
\]
where
\[
 C_0 = C - e_0H_1^{-1}
\]
is again non-negative bur smaller than \(C\), in particular,
\[
\min\sigma(C_0) = 0.
\]
Indeed,
\[
\frac{(C_0\psi,\psi)}{(\psi,\psi)} = \frac{\alpha_0(\phi,\phi)}{(\phi,\phi)}
\frac{(\phi,\phi)}{\|H_1^{1/2}\phi\|^2}, \quad \psi = H_1^{1/2}\phi.
\]
where
\begin{equation}\label{infalpha0}
\inf_\phi \frac{\alpha_0(\phi,\phi)}{(\phi,\phi)} = 0,
\end{equation}
\[
\sup_\psi \frac{(\phi,\phi)}{\|H_1^{1/2}\phi\|^2} < \infty.
\]
In this way we can always extract away the trivial
scalar part \(e_0\) of the perturbation
\(\alpha\) (and similarly for a non-positive \(\alpha\)).
In the following theorem we will therefore suppose that
\begin{equation}\label{infalpha}
\inf_\phi \frac{\alpha_0(\phi,\phi)}{(\phi,\phi)} = 0,\
\mbox{ if \(\alpha\) is non-negative},
\end{equation}\begin{equation}\label{supalpha}
\sup_\phi \frac{\alpha_0(\phi,\phi)}{(\phi,\phi)} = 0,\
\mbox{ if \(\alpha\) is non-positive}.
\end{equation}
Then
\begin{eqnarray}\label{minC}
\min\sigma(C) = 0, \mbox{ if } \alpha \mbox{ is non-negative, }
\end{eqnarray}
\begin{eqnarray}\label{maxC}
\max\sigma(C) = 0, \mbox{ if } \alpha \mbox{ is non-positive. }
\end{eqnarray}

\begin{theorem} \label{window}
Let \((\lambda_-\), \(\lambda_+)\), \(H\), \(\alpha\), \(T\), \(C\) be
as in the previuos lemma and let \(\alpha\) satisfy
(\ref{infalpha},\ref{supalpha})
above.
If the interval
\begin{equation}\label{I}
{\cal I} =
\left(\lambda_- + c_+ (a + b|\lambda_-|),\ \lambda_+
+ c_- (a + b|\lambda_+|)\right)
\end{equation}
where
\begin{equation}\label{cpm}
c_- = \min(\sigma(C)) =
\inf_\psi\frac{\alpha(\psi,\psi)}{\|H_1^{1/2}\psi\|^2},\quad
c_+ = \max(\sigma(C)) =
\sup_\psi\frac{\alpha(\psi,\psi)}{\|H_1^{1/2}\psi\|^2},
\end{equation}
is not void then it is contained in \(\rho(T)\).
\end{theorem}
{\bf Proof.} We supose first that the interval \((\lambda_-,\lambda_+)\)
is finite. Then by virtue of (\ref{minC}) or (\ref{maxC})
this interval must contain
\({\cal I}\).

For every \(d \in{\cal I}\) the complementary projections
\[
P_\pm = \frac{1}{2}(\pm \sign(H - d) + 1))
\]
obviously do not depend on \(d\). In the corresponding matrix representation
we have
\[
(H - d) =
\left[\begin{array}{cc}
(H - d)_+    & 0 \\
0       & -(H - d)_-\\
\end{array}\right],
\]
\[
(H - d)H_1^{-1} =
\left[\begin{array}{cc}
(H - d)_+(a + bH_+)^{-1}    & 0 \\
0       & -(H - d)_-(a + bH_-)^{-1}\\
\end{array}\right],
\]
\[
C =
\left[\begin{array}{cc}
C_{11}       &      C_{12}\\
C_{12}^*     &      C_{22}\\
\end{array}\right],
\]
\[
T - d =
H_1^{1/2}
ZH_1^{1/2}
\]
with
\[
Z =
\left[\begin{array}{cc}
A   &  C_{12}\\
C_{12}^* & -B\\
\end{array}\right],
\]
\[
A = (H - d)_+(a + bH_+)^{-1} + C_{11},\quad
B = (H - d)_-(a + bH_-)^{-1} - C_{22}
\]
so that \(Z^{-1} \in {\cal B}({\cal H})\) implies \(d \in \rho(T)\).
By the obvious identity
\[
Z =
\left[\begin{array}{cc}
1   &  0  \\
C_{12}^*A^{-1} &1\\
\end{array}\right]
\left[\begin{array}{cc}
A   &      0                   \\
0   & -B - C_{12}^*A^{-1}C_{12}\\
\end{array}\right]
\left[\begin{array}{cc}
1   &  A^{-1}C_{12}\\
0   & 1\
\end{array}\right]
\]
we see that \(Z^{-1} \in {\cal B}({\cal H})\)
follows, if both operators \(A,B\)
are positive definite, in particular, if both
\[
(H - d)_+H_{1+}^{-1} + c_-\ \mbox{ and }\ (H - d)_-H_{1-}^{-1} - c_+
\]
are positive definite. This, in turn, is equivalent to
\begin{equation}\label{supc-}
1 + c_-\sup_{\lambda \geq \lambda_+}\frac{a + b|\lambda|}{\lambda - d} > 0
\end{equation}
and
\begin{equation}\label{supc+}
1 - c_+\sup_{\lambda \leq \lambda_-}\frac{a + b|\lambda|}{d - \lambda} > 0.
\end{equation}
Noting that (\ref{sup+}) is valid for any possible value of \(\lambda_-\)
we may rewrite (\ref{supc-}) as
\[
\lambda_+ -d + c_-(a + b\lambda_+) > 0\ \mbox{ \& }\ 1 + c_- b > 0.
\]
Here the second inequality is fullfilled by \(0\leq b < 1,\ |c_-| \leq 1\)
whereas the first is implied by \(d \in {\cal I}\). Now for (\ref{supc-}).
If \(\lambda_- > 0\) then by (\ref{prime2}) and (\ref{prime3})
we have
\[
\sup_{\lambda \leq \lambda_-}\frac{a + b|\lambda|}{d - \lambda}
=
\max\left\{\frac{b\lambda_- + a}{d - \lambda_-},\ b\right\}
\]
and (\ref{supc+}) can be written as
\[
1 > c_+\max\left\{\frac{b\lambda_- + a}{d - \lambda_-},\ b\right\}
\]
which is again guaranteed by \(d \in {\cal I}\).
Here, too, the proof is even simpler, if one of \(\lambda_\pm\) is
infinite. Q.E.D.\\

\begin{remark}\label{lemmata} (i) Neither of the above two
theorems appears to be stronger
or weaker than the other --- in spite of the fact that the interval
\({\cal I}\) from Theorem \ref{window0} is smaller than the one
from Theorem \ref{window}. This lack of elegance is due to the
fact that relative bounds are not shift-invariant. (ii)
Both theorems can be understood as upper-semicontinuity spectral bounds.
According to Theorem \ref{window0} a boundary spectral point
\(\lambda\) cannot move further than \(\pm|\lambda|(a + b|\lambda|)\).
Similarly, by Theorem \ref{window} \(\lambda\) can move as far as
\(\lambda + c_\pm(a + b|\lambda|)\). In particular, {\em the spectrum
moves monotonically} even in spectral gaps: for, say, \(\alpha\)
non-negative,
\begin{equation}\label{I+}
{\cal I} =
\left(\lambda_- + c_+ (a + b|\lambda_-|),\ \lambda_+\right).
\end{equation}
(iii) If \(T = H + A\), \(A\) bounded then
\[
{\cal I} =
\left(\lambda_- + \max\sigma(A),\ \lambda_+ + \min\sigma(A)\right).
\]
\end{remark}

{\bf Bounds for the essential spectra.} The proofs of the
preceding two Theorems
have enough of algebraic structure to be transferable
to the Calkin quotient \({\cal C}^*\)algebra
\({\cal B}({\cal H})/{\cal C}({\cal H})\),
where \({\cal C}({\cal H})\) is the ideal of all compact operators.
Using this we will now derive analogous bounds for the essential spectra.

We first list some simple facts which will be used. Let \({\cal A}\) be
a semisimple \({\cal C}^*\)algebra with the identity \(e\).
If \(p \in {\cal A},\ p \neq e,\) be a projection then the subalgebra
\[
{\cal A}_p = \left\{b \in {\cal A}:\ bp = pb = b\right\}
\]
is again semisimple with the unit \(p\). An element
\(b \in {\cal A}\) is invertible in \({\cal A}_p\), if and only if
in \({\cal A}\) its spectrum has zero as an isolated point and the
corresponding projection is \(q = e - p\). If \({\cal A} =
{\cal B}({\cal H})\) then \({\cal A}_p\) is naturally identified with
\({\cal B}(p{\cal H})\). An element \(b \in {\cal A}\) is called positive,
if its spectrum is non-negative. A sum of two positive elements, one
of which is invertible, is itself positive and invertible.
\begin{prop}\label{schur_alg}
Let \(a = a^* \in {\cal A}\) be invertible and let \(p,q \neq 0\)
be the projections belonging to the positive and the negative part of
\(\sigma(a)\), respectively. Let \(b = b^* \in {\cal A}\)
and \(pbp = qbq = 0\). Then \(a + b\) is invertible.
\end{prop}
{\bf Proof.}  The elements \(ap,aq\) are invertible with the
inverses \(a^{(p)},a^{(q)}\) in \({\cal A}_p,{\cal A}_q\), respectively.
Moreover both \(a^{(p)}\) and \(-a^{(q)}\) are positive.
The fundamental identity (the Schur-complement decomposition)
\[
a + b = za_0z^*
\]
with
\[
z = e + qba^{(p)},\quad z^{-1} = e - qba^{(p)},
\]
\[
a_0 = ap + aq - qba^{(p)}bq
\]
is readily verified. Thus, we have to prove the invertibility of
\(a_0\) in \({\cal A}\). Obviously \(\tilde{a} = aq  - qba^{(p)}bq\)
is invertible in \({\cal A}_q\)
(being a sum of negative elements one of which is invertible).
Denoting by \(\tilde{a}^{(q)}\) its inverse in \({\cal A}_q\)
we have
\[
a_0^{-1} = a^{(p)} + \tilde{a}^{(q)}.
\]
Indeed,
\[
(a^{(p)} + \tilde{a}^{(q)})(ap + \tilde{a})
= (ap + \tilde{a})(a^{(p)} + \tilde{a}^{(q)}) = p + q = e.
\]
Q.E.D.\\

We now prove an analog of Theorem \ref{window} for the essential spectrum.

\begin{theorem}\label{window_ess}
Let \((\lambda_-\), \(\lambda_+)\cap \sigma_{ess}(H) = \emptyset\),
\(\lambda_\pm \in \sigma_{ess}(H)\)
and let \(T = H + \alpha\) satisfy Theorem \ref{pseudoF1} as well as
(\ref{infalpha},\ref{supalpha}), respectively.
If the interval
\begin{equation}\label{I_ess}
{\cal I} =
\left(\lambda_- + c_+ (a + b|\lambda_-|),\
 \lambda_+ + c_- (a + b|\lambda_+|)\right)
\end{equation}
with
\begin{equation}\label{cpm_ess}
c_- = \min(\sigma_{ess}(C))\quad
c_+ = \max(\sigma_{ess}(C))
\end{equation}
is not void then \({\cal I}\cap \sigma_{ess}(T) = \emptyset\).
\end{theorem}
{\bf Proof.} Obviously, \(\sigma_{ess}(C) = \sigma(\hat{C})\),
where
\[
\hat{} \quad : \ {\cal B}({\cal H}) \to {\cal A}
\]
is the Calkin homomorphism. Whenever \(C_\zeta\) is invertible in
\({\cal B}({\cal H})\) and in particular for \(\zeta = i\eta\), \(|\eta|\)
large (\ref{T-zeta}) yields
\begin{equation}\label{pseudoresolvent}
r(\zeta) = \hat{(T - \zeta)^{-1}} =
\hat{H_1^{-1/2}}(\hat{D} + \hat{C})^{-1}\hat{H_1^{-1/2}},
\end{equation}
\[
D = (H - \zeta)H_1^{-1}.
\]
By the spectral mapping principle \(r(\zeta)\) is analytically
continued onto the complement of \(\sigma_{ess}(T)\). This
complement
contains all real \(\zeta = d \in (\lambda_-,\lambda_+)\) for which
\((\hat{D} + \hat{C})^{-1} \in {\cal A}\). Obviously
\[
\sigma(\hat{D}) = f(\sigma_{ess}(H)), \quad
f(\lambda) = \frac{\lambda - d}{a + b|\lambda|}
\]
and \(\hat{D}^{-1} \in {\cal A}\). Let
\(p,q \in {\cal A}\) be the projections corresponding to the
positive and the negative part of the spectrum of
\(\hat{D}\), respectively. As in Theorem \ref{window} one proves that
\[
a_p = p\hat{D} + p\hat{C}p,\ -a_q = -q\hat{D} - q\hat{C}q
\]
are positive and invertible in \({\cal A}_p\), \({\cal A}_q\),
respectively. Now apply Proposition \ref{schur_alg}
to
\[
a = a_p + a_q,\ b = p\hat{C}q + q\hat{C}p\]
thus obtaining the invertibility in \({\cal A}\)
of \(a + b = \hat{D} + \hat{C}\).
Q.E.D.\\

In particular, {\em the essential spectrum depends monotonically on
\(\alpha\)}. Of course, if \(C\) is compact then \(c_\pm = 0\)
and we have \((\lambda_-,\lambda_+) \cap \sigma_{ess}(T) = \emptyset\)
as was known from Theorem \ref{relcomp}.\\

There is an essential-spectrum analog of Theorem \ref{window0} as well:
\begin{theorem}\label{window0_ess}
Let \((\lambda_-\), \(\lambda_+)\cap \sigma_{ess}(H) = \emptyset\),
\(\lambda_\pm \in \sigma_{ess}(H)\)
and let \(T = H + \alpha\) satisfy Theorem \ref{pseudoF1}.
If the interval
\begin{equation}\label{I0_ess}
{\cal I} =
\left(\lambda_- + (a + b|\lambda_-|),\ \lambda_+ -(a + b|\lambda_+|)\right)
\end{equation}
is not void then \({\cal I}\cap \sigma_{ess}(T) = \emptyset\).
\end{theorem}
The proof is similar as above and is omitted.\\

\section{Finite eigenvalues}
All forms in this section will be symmetric.
The following theorem is a necessary tool from the analytic
perturbations which will be repeatedly used later on.
\begin{theorem} \label{holomorphic}
Let \(H\), \(\alpha = \alpha_\epsilon\) for \(\epsilon\)
from an open interval \({\cal I}\) satisfy the conditions of
Theorem \ref{pseudoF} and such that \(\alpha_\epsilon\) is symmetric and
\(C = C_\epsilon\) from (\ref{CC}) is real analytic in
\({\epsilon \in \cal I}\) and
\[
C_{\zeta,\epsilon} = (H - \zeta)H_1^{-1} + C_\epsilon
\]
is invertible in \({\cal B}({\cal H})\) for all \(\zeta\)
from an open set \({\cal O} \subseteq {\cal C}\) and all
\({\epsilon \in \cal I}\).
Then the operator
family \(T_\epsilon = T + \alpha_\epsilon\) is holomorphic in
the sense of \cite{kato},
Ch.~VII, 1. Moreover, the derivative of an isolated
holomorphic eigenvalue
\(\lambda(\epsilon)\) of \(T_\epsilon\) with finite multiplicity
is given by
\begin{equation}\label{lambdaprime}
\lambda'(\epsilon) = \frac{1}{m}
Tr\left((H_1^{1/2}P_\epsilon)^*C'_\epsilon H_1^{1/2}P_\epsilon\right).
\end{equation}
Here \(m,P_\epsilon\) denotes the multiplicity and the spectral
projection on the (total) eigenspace for \(\lambda(\epsilon)\),
respectively.
\end{theorem}
{\bf Proof.} The formula (\ref{lambdaprime}) is plausible being
akin to known analogous expressions from the analytic
perturbation theory (\cite{kato}, Ch.~VII). For completeness
we provide a proof in this more general situation.\footnote{Our case
is close to the holomorphic family of type (C) from
\cite{kato}, Ch.~VII, \S 5.1 where no such details are elaborated.}

Let \(\epsilon_0 \in {\cal I}\) and let \(\Gamma\) be a closed
Jordan curve separating \(\lambda(\epsilon_0)\) from the rest
of \(\sigma(T_{\epsilon_0})\). Let \(\Gamma_1 \subseteq 
\rho(T_{\epsilon_0})\) be another curve connecting \({\cal O}\) and
\(\Gamma\). Take any connected neighbourhood \({\cal O}_0\) of
\(\Gamma\cup\Gamma_1\) with \({\cal O}_0\subseteq 
\rho(T_{\epsilon_0})\).
According to \cite{kato}, Ch.~VII Th.~1.7 
there exists a complex neighbourhood
\({\cal U}_0\) of \(\epsilon_0\) such that 
\((\lambda - T_\epsilon)^{-1}\) is holomorhic in
\({\cal O}_0 \times {\cal U}_0\).

For
\({\lambda \in \cal O}\) and \({\epsilon \in \cal U}_0\) we have
\[
(\lambda - T_\epsilon)^{-1} =
- H_1^{-1/2}C_{\lambda,\epsilon}^{-1}H_1^{-1/2},
\]
\[
\frac{\partial}{\partial \epsilon}(\lambda - T_\epsilon)^{-1}
=
H_1^{-1/2}C_{\lambda,\epsilon}^{-1}
C_\epsilon'
C_{\lambda,\epsilon}^{-1}H_1^{-1/2}.
\]
Note that \({\cal R}(P_\epsilon) \subseteq {\cal Q}\) and hence
\[
H_1^{-1/2}C_{\lambda,\epsilon}^{-1}H_1^{-1/2}
P_\epsilon = \frac{1}{\lambda - \lambda(\epsilon)}P_\epsilon.
\]
By \(H_1^{1/2}P_\epsilon \in {\cal B}({\cal H})\) we have
\[
C_{\lambda,\epsilon}^{-1}H_1^{-1/2}
P_\epsilon =
\frac{1}{\lambda - \lambda(\epsilon)}H_1^{1/2}P_\epsilon,
\]
\begin{equation}\label{PlamPderivative}
P_\epsilon\frac{\partial}{\partial\epsilon}(\lambda - T_\epsilon)^{-1}P_\epsilon =
\frac{1}{(\lambda - \lambda(\epsilon))^2}
(H_1^{1/2}P_\epsilon)^*C_\epsilon' H_1^{1/2}P_\epsilon.
\end{equation}
On the other hand (see \cite{kato})
\[
P_\epsilon = \frac{1}{2\pi i}\int_\Gamma 
(\lambda - T_\epsilon)^{-1}d\lambda,
\]
\[
T_\epsilon P_\epsilon =
\frac{1}{2\pi i}\int_\Gamma \lambda(\lambda - T_\epsilon)^{-1}d\lambda
= \frac{1}{2\pi i}\int_\Gamma 
\lambda P_\epsilon(\lambda - T_\epsilon)^{-1}P_\epsilon d\lambda,
\]
\[
\lambda(\epsilon) = Tr(T_\epsilon P_\epsilon).
\]
Using \(P_\epsilon^2 = P_\epsilon\) and \(P_\epsilon'P_\epsilon = 0\)
we have
\begin{equation}\label{lambdaprev}
\lambda'(\epsilon) = \frac{1}{2\pi im}
Tr\left(\int_\Gamma \lambda P_\epsilon
\frac{\partial}{\partial \epsilon}(\lambda - T_\epsilon)^{-1}
P_\epsilon d\lambda\right)
\end{equation}
(here the integration over \(\lambda\) and the differentiation
over \(\epsilon\) obviously commute).
The formula (\ref{PlamPderivative}) can be analytically continued in
\(\lambda \in {\cal O}_0\)
and inserted into (\ref{lambdaprev}). By using
the obvious identity
\[
\frac{1}{2\pi i}\int_\Gamma \frac{\lambda d\lambda}
{(\lambda - \lambda(\epsilon))^2} = 1
\]
and taking trace we obtain (\ref{lambdaprime}). Q.E.D.\\

The preceding theorem is not general enough to cover
all situations of interest:

\begin{theorem} \label{holomorphic1}
Let \(T_\epsilon = H + \alpha_\epsilon\) be as in Theorem
\ref{holomorphic} above and let \(B(\epsilon)\) be a bounded
symmetric family, analytic in \(\epsilon\). Let, in addition
the set \({\cal O}\) contain a full vertical half-line.
Then the assertions of Theorem
\ref{holomorphic} hold true for
\(T_\epsilon + B(\epsilon)\) but instead of (\ref{lambdaprime})
we have
\begin{equation}\label{lambdaprime1}
\lambda'(\epsilon) = \frac{1}{m}
Tr\left((H_1^{1/2}P_\epsilon)^*C'_\epsilon H_1^{1/2}P_\epsilon
+ P_\epsilon B'(\epsilon)P_\epsilon\right).
\end{equation}
\end{theorem}
{\bf Proof.} We proceed as in the proof of Theorem
\ref{holomorphic} keeping in mind that the formula (\ref{invT-zeta})
is not immediately applicable to \(T_\epsilon + B(\epsilon)\).
We can take \(\Gamma_1\) so as to contain a point
\(\lambda_0 \in {\cal O}\) such that
\(\|B(\epsilon_0)(\lambda - T_{\epsilon_0})^{-1}\| < 1\).
This insures \(\lambda \in \rho(T_\epsilon + B(\epsilon))\)
for \(\lambda \in {\cal O}_1 \subseteq {\cal O}\),
\(\epsilon \in {\cal U}_0\). Then
\[
\frac{\partial}{\partial \epsilon}
(\lambda - T_\epsilon - B(\epsilon))^{-1}
=
\frac{\partial}{\partial \epsilon}\left[
(\lambda - T_\epsilon)^{-1}
\left(1 - B(\epsilon)(\lambda - T_\epsilon)^{-1}\right)^{-1}\right]
\]
\[
=
\frac{\partial}{\partial \epsilon}
(\lambda - T_\epsilon)^{-1}
\left(1 - B(\epsilon)(\lambda - T_\epsilon)^{-1}\right)^{-1} +
\]
\[
(\lambda - T_\epsilon)^{-1}
\left(1 - B(\epsilon)(\lambda - T_\epsilon)^{-1}\right)^{-1}\times
\]
\[
\left(B'(\epsilon)(\lambda - T_\epsilon)^{-1} + 
B(\epsilon)\frac{\partial}{\partial \epsilon}
(\lambda - T_\epsilon)^{-1}\right)\times
\]
\[
\left(1 - (\lambda - T_\epsilon)^{-1}B(\epsilon)\right)^{-1} =
\]
\[
\left(1 - B(\epsilon)(\lambda - T_\epsilon)^{-1}\right)^{-1}
\frac{\partial}{\partial \epsilon}
(\lambda - T_\epsilon)^{-1}
\left(1 - B(\epsilon)(\lambda - T_\epsilon)^{-1}\right)^{-1}
\]
\[
+ \left(\lambda - T_\epsilon - B(\epsilon)\right)^{-1}
B'(\epsilon)
\left(\lambda - T_\epsilon - B(\epsilon)\right)^{-1}
\]
Then using (\ref{PlamPderivative})
\[
P_\epsilon\frac{\partial}{\partial\epsilon}
(\lambda - T_\epsilon - B(\epsilon))^{-1}P_\epsilon =
\frac{1}{(\lambda - \lambda(\epsilon))^2}
\left[
(H_1^{1/2}P_\epsilon)^*C_\epsilon' H_1^{1/2}P_\epsilon
+ P_\epsilon B'(\epsilon)P_\epsilon\right]
\]
which leads to (\ref{lambdaprime1}) as in the theorem above.
Q.E.D.\\

The first application of Theorems \ref{holomorphic},
\ref{holomorphic1} will be
a result on monotonicity. We have to assume that the
spectrum under consideration is sufficiently protected
from unwanted spectral points. We say that a real point \(d\)
is  {\em impenetrable} ({\em essentially impenetrable}) for
a selfadjoint family \(T_\gamma\), \(\gamma\) 
from any set of indices, if \(d \not\in \sigma(T_\gamma)\)
(\(d \not\in \sigma_{ess}(T_\gamma)\)).
\begin{theorem}\label{monotonicity}
Let \(T_\epsilon = H + \alpha_\epsilon\) be analytic
in \(\epsilon \in [\epsilon_0,\epsilon_1]\) in the sense of Theorem 
\ref{holomorphic}.\footnote{
Analyticity in a closed interval
means the same in a complex neighbourhood of that interval.} Let
\(\alpha_\epsilon\) be non-decreasing in \(\epsilon\), let
an open interval \((d,d_1)\) be essentially impenetrable
and one of its ends, say, \(d\) be impenetrable for \(T_\epsilon\).
Let
\[
\lambda_1^1\leq \lambda_2^1\leq \cdots
\]
be the eigenvalues in \((d,d_1)\) of \(T_{\epsilon_1}\)
Then
the spectrum of \(T_{\epsilon_0}\) in \((d,d_1)\)
consists of the eigenvalues which can be ordered as
\[
\lambda_1^0\leq \lambda_2^0\leq \cdots
\]
and they satisfy
\begin{equation}\label{lmonotonicity}
\lambda_k^0 \leq \lambda_k^1,\ k = 1,2,\ldots
\end{equation}
\end{theorem}
{\bf Proof.} For a fixed \(n\) let 
\(\lambda_1^1, \lambda_2^1\leq \cdots \leq \lambda_n^1\)
be the smallest \(n\)
eigenvalues of \(T_{\epsilon_1}\). Any of them
can be analytically continued
to a neighbourhood of \(\epsilon = \epsilon_1\).
By the assumed monotonicity
(we use Theorem \ref{holomorphic}
with \(C_\epsilon' \geq 0\)
as well as the assumed impenetrabilities)
this analytic continuation covers the whole of
\([\epsilon_0,\epsilon_1]\)
i.e.~we obtain analytic non-decreasing functions
\[
d < \lambda_1(\epsilon), \lambda_2(\epsilon),
\ldots \lambda_n(\epsilon) < d_1
\]
as eigenvalues of \(T_\epsilon\). By a permutation, piecewise constant
in \(\epsilon\), we obtain
\[
d < \hat{\lambda}_1(\epsilon) \leq \hat{\lambda}_2(\epsilon)
\leq \cdots \leq \hat{\lambda}_n(\epsilon) < d_1
\]
which are continuous, piecewise analytic,\footnote{
A real function \(f\) is called piecewise analytic on an open
interval \({\cal J}\), if it is real-analytic on
\({\cal J}\setminus {\cal S}\) where \({\cal S}\) is a discrete
set and \(f\) has analytic continuation from each side of
any point from \({\cal S}\), but the two continuatons need not
to coincide.
} still non-decreasing
in \(\epsilon\) and satisfy \(\hat{\lambda}_k(\epsilon_1) = \lambda_k^1\).
By setting \(\epsilon = \epsilon_0\) we obtain \(n\)
eigenvalues of \(T_{\epsilon_0}\)
\[
d < \hat{\lambda}_1(\epsilon_0) \leq \hat{\lambda}_2(\epsilon_0)
\leq \cdots \leq\hat{\lambda}_n(\epsilon_0) < d_1
\]
which obviously satisfy
\[
\hat{\lambda}_k(\epsilon_0) \leq \lambda_k^1,\ k = 1,2,\ldots n.
\]
Then {\em a fortiori}
\begin{equation}\label{fortiori}
\lambda_k^0 \leq \hat{\lambda}_k(\epsilon_0) \leq
\lambda_k^1,\ k = 1,2,\ldots n.
\end{equation}
and \(n\) is arbitrary. The fact that there exists the smallest
eigenvalue \(\lambda_1^0\) of \(T_{\epsilon_0}\) is due to
the impenetrability of the point \(d\). Q.E.D.\\

\begin{remark}\label{monotonicity_rem} Note that, in fact,
the theorem above asserts {\em the existence} of at least
that much eigenvalues of \(T_{\epsilon_0}\) as the \(\lambda_k^1\).
Obviously, if we assume that both interval ends are impenetrable,
then Theorem \ref{monotonicity} applies in both directions
and the eigenvalues of \(T_{\epsilon_0}\)
and \(T_{\epsilon_1}\) have the same
cardinality which is finite.
\end{remark}
\begin{remark}\label{monotonicity1} Theorem \ref{monotonicity}
also holds under the conditions of Theorem \ref{holomorphic1},
if we assume that the form \(\alpha_\epsilon + B(\epsilon)\)
is non-decreasing.
\end{remark}

\begin{cor}\label{Cmonotonicity} Let in Theorem
\ref{monotonicity} \(\alpha_\epsilon = \alpha_0 + \epsilon\alpha_1\)
and let in (\ref{lmonotonicity}) the equality hold for
some \(k\). Then there is \(\psi \neq 0\) with
\begin{equation}\label{alpha1_right}
T_{\epsilon_1}\psi = T_{\epsilon_0}\psi = \lambda_k^1\psi,\quad
\end{equation}
\end{cor}
{\bf Proof.} By the assumption, and using the inequalities
(\ref{fortiori}) from the proof of Theorem
\ref{monotonicity} we obtain
\[
\lambda_k = \hat{\lambda}_k(\epsilon) = \hat{\lambda}_k(\epsilon_0)
= \hat{\lambda}_k(\epsilon_1),\  \mbox{ for all }
 \epsilon \in [\epsilon_0,\epsilon_1].
\]
Thus, \(\hat{\lambda}_k(\epsilon)\) is constant in
\(\epsilon \in [\epsilon_0,\epsilon_1]\). Now
(\ref{lambdaprime}) yields
\begin{equation}\label{Tr=0}
Tr (H_1^{1/2}P_\epsilon)^*C_1H_1^{1/2}P_\epsilon = 0
\end{equation}
for \(\epsilon\) from a neighbourhood of \(\epsilon_1\),
where \(P_\epsilon\) is the (total) projection belonging to
the spectral point \(\lambda_k(\epsilon)\),
\(\epsilon < \epsilon_1\) and
\[
C_\epsilon = C_0 + \epsilon C_1
\]
with \(C_\epsilon' = C_1\) non-negative. Thus, 
(\ref{Tr=0}) implies
\[
C_1H_1^{1/2}P_\epsilon = 0
\]
and, in particular,
\[
\alpha_1(\psi,\phi) = 0,\ \mbox{ for all } \phi \in {\cal Q},
\]
where \(T_{\epsilon}\psi = \lambda_k^1\psi\) for all
\(\epsilon\), in particular,
\(T_{\epsilon_0}\psi = T_{\epsilon_1}\psi = \lambda_k^1\psi
= \lambda_k^0\psi\). Q.E.D.\\

%
The existence of an impenetrable point \(d\) was crucial
in Theorem \ref{monotonicity}. It can be guaranteed by one
of the spectral inclusions, contained in Theorems
\ref{diagonalb}, \ref{window0}, \ref{window}; each of them
contains some restrictions on the size of \(\alpha\) in
comparison to \(H\). Deeper reaching criteria will compare
an 'unknown' \(\alpha\) with a known \(\alpha_0\), which has
desired properties:

\begin{definition}\label{regular}
let \(H\), \(\alpha = \alpha_0 \leq 0\), \({\cal Q}\) be as in
(\ref{alphabound}), (\ref{H1}), (\ref{Q}).\footnote{
Of course, \(\alpha_0 \geq 0\) would do as well. Our
definition of the regularity is, in fact, modeled after 
a standard situation in the applications: the Dirac operator with
the attractive Coulomb potential.} Set
\begin{equation}\label{calA}
{\cal A} = \{\alpha:\ {\cal D}(\alpha) \supseteq {\cal Q},\
|\alpha| \leq c\alpha_0,\ c < 1\}.
\end{equation}
We call \(\alpha_0\) \(H\)-{\em regular}, if the following four
conditions are fulfilled:
\begin{enumerate}
\item Each \(\alpha \in {\cal A}\) satisfies the conditions
of Theorem \ref{pseudoF},
\item \(\sigma_{ess}(H + \alpha) \subseteq \sigma_{ess}(H) =
(-\infty,-m]\cup[m,\infty)\) for some \(m > 0\) and all
\(\alpha \in {\cal A}\),
\item
For some \(\delta > 0\) and all \(\eta\) with 
\(0 \leq \eta < 1\)
\begin{equation}\label{mdelta}
(-m,-m + \delta] \subseteq \rho(H + \eta\alpha_0)
\end{equation}
\item \(\max\sigma(C_0) = 0\), where \(C_0,\ C\) are generated
by (\ref{C}) and \(\alpha_0,\ \alpha\), respectively.
\end{enumerate}
\end{definition}

\begin{theorem}\label{comparison}
Let \(\alpha_0\) be \(H\)-regular and \(\alpha \in {\cal A}\),
\(\alpha \leq 0\).
Then
\[
(-m,-m + \delta] \subseteq \rho(H + \eta\alpha),\quad
0 \leq \eta \leq 1
\]
with \(m,\delta\) from Definition \ref{regular}.
\end{theorem}
{\bf Proof.} Take \(\eta_0 \in (0,1]\)
such that for 
\(a,b\) from (\ref{alphabound}) we have \(\eta_0 b < 1\) 
and
\[
-m + \delta < m + \eta_0c_-^0(a + bm)
\]
where \(c_-^0 = \min\sigma(C_0)\).
Then the conditions of Theorems \ref{pseudoF1}, \ref{window}
hold for \(H + \alpha_{\epsilon,\eta}\) with
\[
\alpha_{\epsilon,\eta}
= (1 - \epsilon)\eta c\alpha_0 + \epsilon \eta \alpha,
\]
(\(c\)  from Def.~\ref{regular}) 
uniformly in \(0 \leq \eta \leq \eta_0\),
\(0 \leq \epsilon \leq 1\); this follows from
\[
c\alpha_0 \leq \alpha_{\epsilon,\eta} \leq 0.
\]
Obviously, \(\alpha_{\epsilon,\eta}\)
belongs to \({\cal A}\) and is non-decreasing in \(\epsilon\) and
non-increasing in \(\eta\). By Theorem \ref{window} we have
\[
(-m, -m + \delta]\subseteq
(-m, m + \eta c_-^0(a + bm)) \subseteq
\rho(H + \alpha_{\epsilon,\eta}),
\]
\(0 \leq \eta \leq \eta_0\), \(0 \leq \epsilon \leq 1\). By
\[
H + \alpha_{\epsilon,\eta} - \zeta
= H_1^{1/2}((H - \zeta)H_1^{-1} + (1 - \epsilon)\eta C
+ \epsilon\eta cC_0)H_1^{1/2}
\]
we see that \(H + \alpha_{\epsilon,\eta}\) is continuous
in the sense of the uniform resolvent topology jointly
in  \(\epsilon,\eta \in [0,1]\) and the same is true for
\(\sigma(H + \alpha_{\epsilon,\eta})\) (see \cite{kato}
Ch.~V.~Th.~4.10). Thus, the set
\[
{\cal S} =\{\eta \in[0,1]:\ (-m,-m + \delta]
\subseteq \rho(H + \alpha_{\epsilon,\eta})\ \mbox{ for all }
\epsilon \in [0,1]\}
\]
is open in \([0,1]\) and it obviously contains \([0,\eta_0]\).
We will prove that the component of \({\cal S}\) containing
\([0,\eta_0]\) is equal to \([0,1]\). If this were not so
then this component
would read \([0,\eta_1)\), \(\eta_0 \leq \eta_1 < 1\). In this case
there would exist an \(\epsilon_1\) such that
\begin{equation}\label{sigmammd}
\sigma(H + \alpha_{\epsilon_1,\eta_1})\cap(-m,-m + \delta]
\neq \emptyset
\end{equation}
whereas
\begin{equation}\label{mmdeltarho}
(-m,-m + \delta] \subseteq \rho(H + \alpha_{\epsilon,\eta})
\end{equation}
for all \(\eta < \eta_1\) and all \(\epsilon \in [0,1]\).

Now, by the mentioned spectral continuity
we would still have
\((-m, -m + \delta) \subseteq \rho(H + \alpha_{\epsilon,\eta_1})\)
for all \(\epsilon \in [0,1]\), more precisely,
\(-m + \delta = \lambda_1(\epsilon_1,\eta_1)\),
where \(\lambda_1(\epsilon,\eta)\) denotes the lowest eigenvalue
of \(H + \alpha_{\epsilon,\eta}\) in \((-m,m)\).
Thus, Theorem \ref{monotonicity}
is applicable to the family
\[
[0,\epsilon_1] \ni \epsilon \mapsto H + \alpha_{\epsilon,\eta_1}
\]
and by (\ref{mdelta}) we would have
\begin{equation}\label{lambdas_all}
-m + \delta < \lambda_1(0,\eta_1) \leq
\lambda_1(\epsilon_1,\eta_1) \leq -m + \delta
\end{equation}
--- a contradiction. Now take in (\ref{mmdeltarho})
\(\epsilon = 1\) which gives
the statement of our theorem. Q.E.D.\\

The theorem above can be regarded as an abstract analog of a result of
W\"{u}st \cite{wuest}, obtained for the Dirac operator
with the Coulomb interaction \(\alpha_0\).

\begin{cor}\label{regreg} If \(\alpha_0\) is regular then
any non-positive \(\alpha \in {\cal A}\) is regular also.
\end{cor}
\begin{cor}\label{less} Let \(\alpha_0\) be \(H\)-regular
and \(0\geq \beta \geq \alpha \in {\cal A}\). Then the spectrum
of \(H + \alpha\), \(H + \beta\) in \((-m,m)\)
consists of the eigenvalues
\[
\lambda_1 \leq \lambda_2 \leq \cdots
\]
\[
\mu_1 \leq \mu_2 \leq \cdots
\]
respectively, and
\begin{equation}\label{lmless}
\lambda_k \leq \mu_k,\quad k = 1,2\ldots
\end{equation}
holds.
\end{cor}
{\bf Proof.} By Theorem \ref{comparison} we have
\((-m, -m + \delta) \subseteq \rho(T_\epsilon)\)
where
\[
T_\epsilon = H + (1 - \epsilon)\alpha + \epsilon\beta,\quad
0 \leq \epsilon \leq 1,
\]
\[
(1 - \epsilon)\alpha + \epsilon\beta \in {\cal A}.
\]
Now Theorem \ref{monotonicity} applies and (\ref{lmless})
follows. Q.E.D.

\begin{theorem}\label{finebound}
Let \(\alpha_0\) be \(H\)-regular and let
\begin{equation}\label{fineboundalpha}
|\alpha - c\alpha_0| \leq -\epsilon c\alpha_0,
\end{equation}
\[
\epsilon < \min\{1, \frac{1}{c} - 1\}.
\]
Then the spectrum of \(H + \alpha\) in \((-m,m)\)
consists of the eigenvalues
\[
\mu_1 \leq \mu_2 \leq \cdots
\]
and they satisfy
\begin{equation}\label{finelmless}
\lambda_k((1 + \epsilon)c) \leq \mu_k \leq
\lambda_k((1 - \epsilon)c) \quad k = 1,2\ldots
\end{equation}
where 
\[
\lambda_1(\eta) \leq \lambda_2(\eta) \leq \cdots
\]
are the eigenvalues of \(H + \eta\alpha_0\).

\end{theorem}
{\bf Proof.} (\ref{fineboundalpha}) can be written as
\[
(1 + \epsilon) c\alpha_0 \leq \alpha \leq
(1 - \epsilon) c\alpha_0,
\]
from which it follows \(\alpha \leq 0\) and
\[
0 \geq \alpha \geq c_1\alpha_0
\]
with \(c_1 = (1 + \epsilon) c < 1\). Thus, 
\(\alpha \in {\cal A}\). Now apply
Corollary \ref{less} to the operators
\[
H + (1 + \epsilon) c\alpha_0,\ H + \alpha,\
H + (1 - \epsilon) c\alpha_0
\]
and (\ref{finelmless}) follws. Q.E.D.
\begin{remark}\label{sharp} The estimates (\ref{finelmless})
are sharp: by taking the perturbation
\(\beta = (1 \pm \epsilon)\alpha\) the equality on the
respective side in (\ref{finelmless}) is obtained. The
bound (\ref{finelmless}) is particularly useful, if the eigenvalues
\(\lambda(\eta)\) are explicitly known as functions of
\(\eta\) as is the case with several important quantum mechanical
systems.
\end{remark}

Let us now turn to the promised bound (\ref{ev_bound}).
We will combine the monotonicity
from Theorem \ref{monotonicity} with one of the spectral inclusion
results above to insure the necessary impenetrabilities. There are
quite few of the latter, so we will present the most
typical cases.

\begin{theorem}\label{general_bound0}
Let \(H,\alpha,T,C\) be as in Theorem \ref{pseudoF1} and
\(\alpha\) symmetric. Let \({\cal I} =
(\lambda_{--},\lambda_{++})\)
be an essential spectral gap for \(H\) and \(\lambda_{-+}\)
the lowest eigenvalue of \(H\) in \({\cal I}\) such that
the open interval
\begin{equation}\label{sep0-}
{\cal I}_- = (\lambda_{--} + a + b|\lambda_{--}|,\
 \lambda_{-+} - a - b|\lambda_{-+}|)
\end{equation}
is non-void. Furthermore, let either
\begin{itemize}
\item[(i)] \(\lambda_{+-}\) be
the highest eigenvalue of \(H\) in \({\cal I}\) such that
the open interval
\begin{equation}\label{sep0+}
{\cal I}_+ = (\lambda_{+-} + a + b|\lambda_{+-}|,\
\lambda_{++} - a - b|\lambda_{++}|)
\end{equation}
is non-void or
\item[(ii)] the form \(\alpha\) satisfy the conditions
of Theorem \ref{relcomp}.
\end{itemize}
By
\[
\lambda_1 = \lambda_{-+} \leq \lambda_2 \leq \cdots
\]
denote the (finite or infinite) sequence of the eigenvalues
of \(H\) in \({\cal I}\). Set
\[
\tilde{\lambda} =
\left\{\begin{array}{lr}
\lambda_{++} - a - b|\lambda_{++}|,& \mbox{ in case (i)} \\
\lambda_{++},                      & \mbox{ in case (ii)}\\
\end{array}\right.
\]
Then the spectrum of \(T\) in
\(\widetilde{{\cal I}} =
(\lambda_{--} + a + b|\lambda_{-}|,\tilde{\lambda})\)
consists of the eigenvalues
\[
\mu_1 \leq \mu_2 \leq \cdots
\]
and they satisfy
(\ref{ev_bound}) in the following sense: in the case (i)
\(\lambda\)'s and \(\mu\)'s have the same cardinality and
(\ref{ev_bound}) holds for all of them whereas in the case
(ii) (\ref{ev_bound}) holds as long as
\(\lambda_k + a + b|\lambda_k| < \lambda_{++}\).
\end{theorem}
{\bf Proof.} We introduce an auxiliary family
\[
\tilde{T}_\epsilon = H + \tilde{\alpha}_\epsilon,
\]
\[
\tilde{\alpha}_\epsilon = \epsilon (a + b\hat{h}) +
(1 - \epsilon)\alpha,\quad 0 \leq \epsilon \leq 1,
\]
where \(\hat{h}\) is the closed form belonging to the operator
\(|H|\). This family satisfies the conditions
of Theorem \ref{pseudoF1} uniformly in \(\epsilon\):
\[
|\tilde{\alpha}_\epsilon(\psi,\psi)| \leq
\epsilon(a + b\hat{h})(\psi,\psi) +
(1 - \epsilon)(a + b\hat{h})(\psi,\psi)
\]
\[
= (a + b\hat{h})(\psi,\psi)
\]
and the operator \(\tilde{C}_\epsilon\), constructed
from \(\tilde{T}_\epsilon\) according to (\ref{C})
is here given by
\[
(\tilde{C}_\epsilon\psi,\phi) =
\tilde{\alpha}_\epsilon(H_1^{-1/2}\psi,H_1^{-1/2}\phi)
=
((\epsilon + (1 - \epsilon)C)\psi,\phi),
\]
so, \(\tilde{C}_\epsilon = \epsilon + (1 - \epsilon)C\)
is holomorphic with \(\|\tilde{C}_\epsilon\| \leq 1\)
and \(\tilde{C}_\epsilon' = 1 - C\) is non-negative.
In particular, \(\tilde{T}_\epsilon\) fulfills
the conditions of Theorem \ref{holomorphic}
as well as Theorem \ref{window0}, uniformly
in \(\epsilon \in [0,1]\). Thus, \({\cal I}_-\)
is impenetrable for \(\tilde{T}_\epsilon\).
We now show that any open interval \((d,d_1)\)
with \(d \in {\cal I}_-\) and
\[
d_1\ \left\{\begin{array}{rr}
\in {\cal I}_+,  & \mbox{ in case (i)}\\
= \lambda_{++},  & \mbox{ in case (ii)}\\
\end{array}\right.
\]
is essentially impenetrable for \(\tilde{T}_\epsilon\).
To this end we introduce another auxiliary family
\[
H + \epsilon(a + b\hat{h}) = H + \epsilon H_1
\]
to which both Theorem \ref{pseudoF1} and Theorem \ref{window0}
hold, again uniformly in \(\epsilon\). Therefore its spectrum 
in \(\widetilde{{\cal I}}\) consists of the eigenvalues
\[
\lambda_1 + \epsilon(a + b|\lambda_1|)
\leq \lambda_2 + \epsilon(a + b|\lambda_2|)
\leq \cdots
\]
In particular, \(\tilde{{\cal I}}\) is essentially impenetrable
for \(H + \epsilon H_1\).
The form sum \(\tilde{T}_\epsilon = H + \tilde{\alpha}_\epsilon\)
can obviously be represented as another form sum
\begin{equation}\label{tildeTform}
\tilde{T}_\epsilon = (H + \epsilon H_1) + (1 - \epsilon)\alpha
\end{equation}
again in the sense of Theorems \ref{pseudoF} and \ref{pseudoF1}.
Indeed, using the the functional calculus we obtain the operator
inequality (a proof is
provided in the Appendix)
\begin{equation}\label{Hinequality}
H_1 \leq \frac{a}{1 - \epsilon b} + \frac{b|H + \epsilon H_1|}
{1 - \epsilon b},
\end{equation}
hence
\[
|(1 - \epsilon)\alpha(\psi,\psi)| \leq
\frac{1 - \epsilon}{1 - \epsilon b}a\|\psi\|^2 +
\frac{(1 - \epsilon)b\||H + \epsilon H_1|^{1/2}\psi\|^2}
{1 - \epsilon b}
\]
\[
\leq a\|\psi\|^2 + b\||H + \epsilon H_1|^{1/2}\psi\|^2.
\]
Furthermore, the form sum (\ref{tildeTform}) satisfies
the conditions of Theorem \ref{relcomp}.
As a matter of fact, the operator \(\tilde{C}\), defined by
\[
(\tilde{C}\psi\phi) = \alpha((a + b|H + \epsilon H_1|)^{-1/2}\psi,
(a + b|H + \epsilon H_1|)^{-1/2}\phi)
\]
satisfies
\[
(a + b|H + \epsilon H_1|)^{-1}\tilde{C} = BH_1^{-1}CB
\]
where \(H_1^{-1}C\) is known to be compact and by (\ref{Hinequality})
\[
B = H_1^{1/2}(a + b|H + \epsilon H_1|)^{-1/2}
\]
is bounded. Thus, \((d,d_1)\) is essentially impenetrable
for \(\tilde{T}_\epsilon\) in the case (ii).
The case (i) is even simpler: due
to the impenetrability from both sides for \(H + \epsilon H_1\)
the cardinalities of the eigenvalues of \(H\) and \(H + H_1\)
are finite and equal, the same is then true of \(T\) and
\(H + H_1\) now  due to the impenetrability from both
sides for \(\tilde{T}_\epsilon\).
Now all conditions of Theorem \ref{monotonicity} are fulfilled
for the family \(\tilde{T}_\epsilon\) for which
\(\tilde{T}_0 = T\) and \(\tilde{T}_1 = H + H_1\).
Hence the eigenvalues
of \(T\)
in \(\tilde{{\cal I}}\) are
\[
\mu_1 \leq \mu_2 \leq \cdots,
\]
they are at least as much as those
\(\lambda_k + a + b|\lambda_k|\) which are smaller than
\(\tilde{\lambda}\) and they satisfy the right hand side of
(\ref{ev_bound}). To obtain the other
we use the form \(\tilde{\alpha}_\epsilon
= -\epsilon (a + b\hat{h}) +
(1 - \epsilon)\alpha\). Q.E.D.\\

\begin{remark}\label{rem:general}(i)
In the proof above the right hand side of the inequality
(\ref{ev_bound}) had to be proved first because this step guarantees
{\em the existence} of the perturbed eigenvalues. This
asymetry is natural and is due to the fact that in general
only the left end of the 'window' \((d,d_1)\) is assumed
as impenetrable (case (ii)). The other direction is
handled by considering \(H = - H\). (ii) The restrictive condition
that \(\lambda_k + a + b|\lambda_k|\) be smaller than
\(\tilde{\lambda}\) is trivially fulfilled for all \(k\),
if \(\lambda_{++} = \infty\).
\end{remark}

An analogous result holds under the conditions of Theorem
\ref{window}.
\begin{theorem}\label{general_bound}
Let \(H,\alpha,T,C\) be as in Theorem \ref{pseudoF1} and
\(\alpha\) symmetric and let, in
addition, \(\alpha\) satisfy (\ref{infalpha}), (\ref{supalpha}) with
\(c_\pm\) from (\ref{cpm}). Let \({\cal I} =
(\lambda_{--},\lambda_{++})\)
be an essential spectral gap for \(H\) and \(\lambda_{-+}\)
the lowest eigenvalue of \(H\) in \({\cal I}\) such that
the open interval
\begin{equation}\label{sep-}
{\cal I}_- = (\lambda_{--} + c_+(a + b|\lambda_{--}|),\
 \lambda_{-+} + c_-(a + b|\lambda_{-+}|))
\end{equation}
is non-void. Furthermore, let either
\begin{itemize}
\item[(i)] \(\lambda_{+-}\) be
the highest eigenvalue of \(H\) in \({\cal I}\) such that
the open interval
\begin{equation}\label{sep+}
{\cal I}_+ = (\lambda_{+-} + c_+(a + b|\lambda_{+-}|),\
 \lambda_{++} + c_-( a + b|\lambda_{++}|))
\end{equation}
is non-void or
\item[(ii)] the form \(\alpha\) satisfy the conditions
of Theorem \ref{relcomp}.
\end{itemize}
By
\[
\lambda_1 = \lambda_{-+} \leq \lambda_2 \leq \cdots
\]
denote the (finite or infinite) sequence of the eigenvalues
of \(H\) in \({\cal I}\). Set
\[
\tilde{\lambda} =
\left\{\begin{array}{lr}
\lambda_{++} + c_-( a + b|\lambda_k|),& \mbox{ in case (i)} \\
\lambda_{++},                         & \mbox{ in case (ii)}\\
\end{array}\right.
\]
Then the spectrum of \(T\) in
\(\widetilde{{\cal I}} =
(\lambda_{--} + a + b|\lambda_{-}|,\tilde{\lambda})\)
consists of the eigenvalues
\[
\mu_1 \leq \mu_2 \leq \cdots
\]
and they satisfy
\begin{equation}\label{ev_bound+-}
\lambda_k + c_-( a + b|\lambda_k|)
\leq \mu_k \leq \lambda_k + c_+( a + b|\lambda_k|).
\end{equation}
in the following sense: in the case (i)
\(\lambda\)'s and \(\mu\)'s have the same cardinality and
(\ref{ev_bound}) holds for all of them whereas in the case
(ii) (\ref{ev_bound}) holds as long as
\(\lambda_k + c_+( a + b|\lambda_k|) < \lambda_{++}\).
\end{theorem}

We omit the proof, it follows the lines of the one
of Theorem \ref{general_bound0} above.
The only difference in the proof is the form \(\tilde{\alpha}_\epsilon\)
which now reads
\[
\tilde{\alpha}_\epsilon = c_\pm\epsilon (a + b\hat{h}) +
(1 - \epsilon)\alpha.
\]
Also, Remark
\ref{rem:general} applies accordingly.
\begin{remark}\label{alphapositive}
If in the preceding theorem the form \(\alpha\) is
non-negative then the bound (\ref{ev_bound+-}) reads
\begin{equation}\label{ev_bound2_positive}
0 \leq \mu_k - \lambda_k \leq c_+(a + b|\lambda_k|)
\leq a + b|\lambda_k|.
\end{equation}
\end{remark}

The preceding theorems cover perturbation estimates already known:
by setting \(a = 0\) the bound (\ref{ev_bound}) was obtained in
\cite{vessla1}
for finite matrices. Also by setting \(b = 0\) we have
\(T = H + A\), \(A \in {\cal B}({\cal H})\), \( C = A/a\); here
(\ref{ev_bound+-}) gives the mentioned bound
(\ref{ev_bound_A+-}).
Both (\ref{ev_bound}) and (\ref{ev_bound+-}) are
sharp, they obviously become equalities on scalars.

Positioning of an impenetrable point is user dependent; usually a
most convenient choice is to take broad
spectral gaps. In the most notorious case of a
positive definite \(H\) with
a compact inverse the impenetrability from below
is trivially fulfilled.\\

The proofs of Theorem \ref{general_bound0}
and \ref{general_bound} consist of two main ingredients:
\begin{enumerate}
\item upper semicontinuity bounds for general spectra from Theorem \ref{window0_ess}, \ref{window_ess} and
\item lower semicontinuity bounds for finite eigenvalues, obtained by
the construction of monotone holomorphic operator families.
\end{enumerate}

So, we may say that in order to fully control the eigenvalues in a
gap by using (\ref{ev_bound}) or (\ref{ev_bound+-})
have to 'pay a price', that is, the perturbation should be so small
as to insure that the impenetrability conditions (\ref{sep0-}),
(\ref{sep0+}), (\ref{sep-}), (\ref{sep+}), respectively, are
fulfilled. These expressions as well as the estimates in
(\ref{ev_bound}) or (\ref{ev_bound+-}) use the same bound \(\pm( a +
b|\lambda|)\), so the price is completely adequate. This fact may be
seen as a mark of the naturality of the obtained bounds.\\

\section{Appendix}
{\bf Proof of (\ref{maxf}).}
Obviously the point \(\xi = 0\) is a local minimum of 
\(\psi(\cdot,a,b,\lambda,\eta)\). By
\[
\psi(-\xi,a,b,\lambda,\eta) = \psi(\xi,a,b,-\lambda,\eta)
\]
it is sufficient to take \(\lambda \geq 0\). We distinguish two
cases.

\(\xi \geq 0\):

\[
\psi_\xi = \frac{-\xi(a + \lambda b) + (\lambda^2 + \eta^2)b + \lambda a}
{((\xi - \lambda)^2 + \eta^2)^{3/2}}.
\]
The maximum is reached at
\[
\xi = \xi_0 = \lambda + \frac{\eta^2 b}{a + \lambda b}
\]
and it is equal to
\[
\psi(\xi_0,a,b,\lambda,\eta) = \frac{1}{|\eta|}
\sqrt{(a + \lambda b)^2 + \eta^2 b^2}
\]
and this is (\ref{maxf}).\\

\(\xi \leq 0\):

\[
\psi_\xi = \frac{-\xi(a - \lambda b) - (\lambda^2 + \eta^2)b + \lambda a}
{((\xi - \lambda)^2 + \eta^2)^{3/2}}.
\]
The maximum is reached at
\[
\xi = \xi_1 = \lambda - \frac{\eta^2 b}{a - \lambda b}
\]
and it is equal to
\[
\psi(\xi_1,a,b,\lambda,\eta) = \frac{1}{|\eta|}
\sqrt{(a - \lambda b)^2 + \eta^2 b^2},
\]
provided that \(a > \lambda b\) and
\[
\lambda \leq \frac{b\eta^2}{a - \lambda b},
\]
otherwise the maximum is reached on the boundary \(\{-\infty,0\}\).
All three values are obviously less than (\ref{maxf})
which is the sought global maximum.\\

{\bf Proof of (\ref{Hinequality}).}
For real \(\lambda\) we have
\[
|\lambda + \epsilon (a + b|\lambda|)| =
\left\{\begin{array}{rr}
\lambda + \epsilon (a + b\lambda,    & \lambda \geq 0 \\
|\lambda + \epsilon (a - b\lambda)|, & \lambda \leq 0 \\
\end{array}\right.
\]
Thus, for \(\lambda \geq 0\)
\[
\lambda =
\frac{-\epsilon a}{1 + b\epsilon} +
\frac{|\lambda + \epsilon (a + b|\lambda|)|}{1 + b\epsilon} \leq
\]
\[
\frac{\epsilon a}{1 - b\epsilon} +
\frac{|\lambda + \epsilon (a + b|\lambda|)|}{1 - b\epsilon}
\]
and for \(\lambda \leq 0\)
\[
|\lambda + \epsilon (a + b|\lambda|)| =
|\epsilon a + \lambda(1 - b\epsilon)| \geq
-\epsilon a - \lambda(1 - b\epsilon)
\]
hence
\[
-\lambda \leq
\frac{\epsilon a}{1 - b\epsilon} +
\frac{|\lambda + \epsilon (a + b|\lambda|)|}{1 - b\epsilon}.
\]
Altogether
\[
|\lambda| \leq
\frac{\epsilon a}{1 - b\epsilon} +
\frac{|\lambda + \epsilon (a + b|\lambda|)|}{1 - b\epsilon}.
\]
Taking corresponding functions of \(H\) we obtain
 (\ref{Hinequality}).


\begin{thebibliography} {99}

\bibitem{dolbeault} Dolbeault, J., Esteban, M.~J., S\`{e}r\`{e}, E.,
Variational characterisation for eigenvalues of Dirac operators,
Preprint mp-arc 98-177, to appear in Calc.~Var.~and PDE.

\bibitem{faris} Faris, W.G,. Self-adjoint operators.
Lecture Notes in Mathematics, Vol. 433. Springer-Verlag,
 Berlin-New York, 1975.
\bibitem{kato}Kato, T., Perturbation Theory for Linear Operators,
Springer, Berlin 1966.
\bibitem{griesemer} Griesemer, M., Lewis, R.T., Siedentop, H.,
A minimax principle for eigenvalues in spectral gaps:
Dirac operators with Coulomb
potentials, Documenta Mathematica {\bf 4} (1999) 275-283.
\bibitem{gkmv} Grubi\v{s}i\'{c}, L., Kostrykin, V., Makarov, K.A.,
Veseli\'{c}, K., On the perturbation theory for quadratic forms,
in preparation.

\bibitem{nenciu} Nenciu, G., Self-adjointness and invariance
of the essential spectrum for Dirac operators defined as
quadratic forms, Comm.~Math.~Phys.~{\bf 48} (1976) 235--247.

\bibitem{thaller} Thaller, B. The Dirac Equation, Springer 1992.
\bibitem{vesfact} Veseli\'{c}, K., Perturbation theory for the eigenvalues
of factorised symmetric matrices, LAA {\bf 309} (2000) 85-102.


\bibitem{vessla1} Veseli\'{c}, K., Slapni\v{c}ar, I., Floating point
 perturbations of Hermitian
matrices, Linear Algebra Appl. {\bf 195} (1993) 81-116.
\bibitem{winklmeier}
Winklmeier, M., The Angular Part of the Dirac Equation
in the Kerr-Newman Metric:
Estimates for the Eigenvalues, Ph.~D.~Thesis 2005.
\bibitem{wuest}
W\"{u}st, R.,
Dirac operations with strongly singular potentials. Distinguished self-adjoint extensions constructed with a spectral gap theorem and cut-off potentials.
Math. Z. {\bf 152} (1977), no. 3, 259--271. 
\end{thebibliography}
\end{document}